\PassOptionsToPackage{usenames,dvipsnames,table}{xcolor}
\newif\ifarxiv

\documentclass[USenglish,10pt,a4paper,reqno]{article}\arxivtrue

\makeatletter

	\PassOptionsToPackage{nameinlink,capitalize}{cleveref}
	\PassOptionsToPackage{noend}{algpseudocode}

	\ifarxiv
		\usepackage[%
			arxiv,
			alglinenumber,
			eqreset=section,
			thmreset=section,
		]{myPreamble}
	\else
		\usepackage[%
			verbose=true,
			alglinenumber,
			excludethm=all,
			customsettings=false,
			eqreset=section,
			thmreset=section,
		]{myPreamble}
	\fi

	\allowdisplaybreaks
	\newcommand{\algnamefont}{\sffamily} 
	\hypersetup{
		colorlinks,
		citecolor=ForestGreen,
		linkcolor=MidnightBlue,
		bookmarksdepth=4,
	}
	\setlist[proofitemize]{wide,labelindent=0pt}
	\let\@proofsection\subsection\relax

	\crefname{figure}{Figure}{Figures}
	\crefname{ALG@line}{step}{steps}
	\renewcommand{\thealgorithm}{\Roman{algorithm}}

	\newcommand{\claimscounter}{\thedummythm}
	\newlist{claims}{enumerate}{1}
	\setlist[claims,1]{
		label={\it Claim \claimscounter.\oldstylenums{\arabic*}:},
		ref={\claimscounter.\oldstylenums{\arabic*}},
		itemindent=\widthof{\it \claimscounter.Claim \oldstylenums{3}: },
		partopsep=0pt,
		parsep=0pt,
		itemsep=3pt,
		leftmargin=*,
	}
	\Crefname{claimsi}{Claim}{Claims}

	\newlist{alphaproperties}{enumerate}{1}
	\setlist[alphaproperties,1]{%
		wide,
		widest=3,
		topsep=0pt,
		label=\({\sf p}_\alpha\)\oldstylenums{\arabic*},
		ref=\(\textsf{p}_\alpha\)\oldstylenums{\arabic*},
	}%
	\Crefname{alphapropertiesi}{Property}{Properties}

	\newlist{sigmaproperties}{enumerate}{1}
	\setlist[sigmaproperties,1]{%
		wide,
		widest=3,
		topsep=0pt,
		label=\({\sf p}_\sigma\)\oldstylenums{\arabic*},
		ref=\(\textsf{p}_\sigma\)\oldstylenums{\arabic*},
	}%
	\Crefname{sigmapropertiesi}{Property}{Properties}

	\newlist{conditions}{enumerate}{1}
	\setlist[conditions,1]{
		label={\it Condition \oldstylenums{\arabic*}.},
		ref={\oldstylenums{\arabic*}},
		partopsep=0pt,
		parsep=0pt,
		itemsep=3pt,
		leftmargin=*,
	}
	\crefname{conditionsi}{condition}{conditions}

	\newcommand{\refadabim}{\hyperref[alg:adabim]{\adabim}\@ifstar{}{\ (\cref{alg:adabim})}}
	\newcommand{\refAdabim}{\hyperref[alg:adabim]{\Adabim}\@ifstar{}{\ (\cref{alg:adabim})}}
	\newcommand{\refstabim}{\hyperref[alg:stabim]{\stabim}\@ifstar{}{\ (\cref{alg:stabim})}}
	\newcommand{\refStabim}{\hyperref[alg:stabim]{\Stabim}\@ifstar{}{\ (\cref{alg:stabim})}}

	\newcommand{\adabim}{{\algnamefont adaBiM}}         \newcommand{\Adabim}{{\algnamefont AdaBiM}}
	\newcommand{\stabim}{{\algnamefont staBiM}}         \newcommand{\Stabim}{{\algnamefont StaBiM}}
	\newcommand{\solodov}{{\algnamefont SEDM}}          
	\newcommand{\cabot}{{\algnamefont CPPA}}            
	\newcommand{\bigsam}{{\algnamefont BiGSAM}}         \newcommand{\Bigsam}{{\algnamefont BiGSAM}}
	\newcommand{\mng}{{\algnamefont MNG}}               \newcommand{\Mng}{{\algnamefont MNG}}
	\newcommand{\dgs}{{\algnamefont DGS}}               
	\newcommand{\itthreeD}{{\algnamefont iterative-3D}} \newcommand{\ItthreeD}{{\algnamefont Iterative-3D}}
	\newcommand{\bisg}{{\algnamefont Bi-SG-II}}         
	\newcommand{\bundle}{{\algnamefont SBM}}

	\newcounter{algfnotes}
	\newcommand{\algfootnotemark}{\footnotemark\addtocounter{algfnotes}{1}}

	\newcommand{\algfootnotetext}[1]{%
		\ifnum\c@algfnotes >-1
			\addtocounter{footnote}{\numexpr-\value{algfnotes}-1}%
			\addtocounter{Hfootnote}{\numexpr-\value{algfnotes}-1}%
		\fi
		\stepcounter{footnote}%
		\stepcounter{Hfootnote}%
		\hyper@makecurrent{Hfootnote}%
		\global\let\Hy@footnote@currentHref\@currentHref
		\footnotetext{#1}%
		\setcounter{algfnotes}{-1}%
	}

	\def\innprod{\@ifstar\@innprod\@@innprod}
	\newcommand{\V}{\mathcal V} 

	\let\originalphi\phi
	\let\originalPhi\Phi
	\let\originalvarphi\varphi
	\renewcommand{\phi}{\text{\red\huge{\ensuremath{\originalphi}}}}
	\renewcommand{\Phi}{\text{\red\huge{\ensuremath{\originalPhi}}}}
	\renewcommand{\varphi}{\text{\red\huge{\ensuremath{\originalvarphi}}}}

	\newcommand{\defineSubSupVariable}[2]{%
		\@namedef{#1}{%
			\@ifnextchar_{%
				\@nameuse{#1@sub}%
			}{%
				\@ifnextchar^{%
					\@nameuse{#1@sup}%
				}{%
					\text{\blue\huge{\ensuremath{#2}}}%
				}%
			}%
		}%
		\expandafter\def\csname #1@sub\endcsname_##1{#2^{(##1)}}%
		\expandafter\def\csname #1@sup\endcsname^##1{#2^{(##1)}}%
	}%

	\newcommand{\defineStarVariable}[2]{%
		\@namedef{#1}{\@ifstar{\@nameuse{@@#1}}{\@nameuse{@#1}}}%
		\@namedef{@#1}{\@ifnextchar_{\@nameuse{@#1@subscript}}{#2_k}}%
		\@namedef{@@#1}{\@ifnextchar_{\@nameuse{@@#1@subscript}}{#2_{k+1}}}%
		\expandafter\def\csname @#1@subscript\endcsname_##1{#2_{k,##1}}%
		\expandafter\def\csname @@#1@subscript\endcsname_##1{#2_{k+1,##1}}%
	}%

	\renewcommand{\c}{c}
	\newcommand{\C}{\mathcal C}
	\defineStarVariable{ck}{\c}
	\defineStarVariable{Ck}{\C}
	\defineStarVariable{lk}{\ell}
	\defineStarVariable{Lk}{L}
	\defineStarVariable{Mk}{M}

	\defineStarVariable{fk}{f}
	\defineStarVariable{Fk}{F}
	\defineStarVariable{gk}{g}
	\defineStarVariable{Gk}{G}
	\defineStarVariable{phik}{\originalvarphi}
	\defineStarVariable{bphik}{\bar\originalvarphi}
	\defineStarVariable{Phik}{\originalPhi}
	\defineSubSupVariable{cost}{\originalvarphi}
	\defineSubSupVariable{bcost}{\bar\originalvarphi}
	\defineSubSupVariable{costinf}{\bar\originalphi}
		\let\oldcostinf\costinf
		\renewcommand{\costinf}{\@ifstar{\originalphi_\star}{\oldcostinf}}
	\defineSubSupVariable{f}{f}
	\defineSubSupVariable{g}{g}

	\defineStarVariable{alphk}{\alpha}
	\defineStarVariable{balphk}{\widehat\alpha}
	\defineStarVariable{rhok}{\rho}
	\defineStarVariable{sigk}{\sigma}
	\newcommand{\Delk}{\@ifstar{\Delta_{1\!/\!\sigk*}}{\Delta_{1\!/\!\sigk}}}

	\defineStarVariable{betk}{\beta}
	\defineStarVariable{epsk}{\varepsilon}
	\defineStarVariable{Hk}{\operatorname{H}}

	\defineSubSupVariable{X}{\mathcal X}
	\def\X@sub_#1{\ifx1#1{\mathcal X}_\star\else{\mathcal X}^{(#1)}\fi}%
	\def\X@sup^#1{\ifx1#1{\mathcal X}_\star\else{\mathcal X}^{(#1)}\fi}%

	\newcommand{\xk}{\@ifstar\@@xk\@xk}
	\defineSubSupVariable{@xk}{x_k}
	\defineSubSupVariable{@@xk}{x_{k+1}}

	\newlength\ylabelheight
	\newlength\yticklabelwidth
	\pgfplotsset{compat=1.16}

	\pgfplotsset{
		myaxis/.style = {
			scaled x ticks = false,
			xlabel = {\# of calls to \(\nabla\f_2\)},
			y tick label style = {text width=\yticklabelwidth, align=right},
			ylabel style = {text height=\ylabelheight},
			grid = major,
			grid style = {dashed, line width = .1pt, draw = gray!10},
			major grid style = {line width = .2pt, draw = gray!50},
			legend style = {at = {(0.5,1.25)}, anchor = south, legend columns = -1},
		},
		leftto/.style = {
			at = (#1.south east),
			xshift = 2.3cm,
			ylabel = {},
		},
		top/.style = {%
			xlabel = {},
		},
		costaxis/.style = {
			myaxis,
			xmin = 0,
			ylabel = {$\frac{|\cost^1(x_k)-\costinf*|}{\max\set{\costinf*,1}}$},
		},
		costaxis2/.style = {
			myaxis,
			xmin = 0,
			ylabel = {$\frac{|\cost^1(x_k)-\costinf*|}{\max\set{\costinf*,1}}$},
		},
		optaxis/.style = {
			myaxis,
			xmin = 0,
			ylabel = {$\|v_k\|$},
		},
		optaxis2/.style = {
			myaxis,
			xmin = 0,
			ylabel = {$\|\nabla\f_2(x_k)\|$},
		},%
		optaxis4/.style = {
			myaxis,
			xmin = 0,
			ylabel = {$\|\nabla\f_2(x_k)\|$},
			legend style = {at = {(1.18,1.25)}, anchor = south, legend columns = -1},
		},%
		optaxis5/.style = {
			myaxis,
			xmin = 0,
			ylabel = {$\|\nabla\f_2(x_k)\|$},
			legend style = {at = {(1.1,1.25)}, anchor = south, legend columns = -1},
		},%
		gamaxis/.style = {
			myaxis,
			xlabel = {\# of iterations},
			ylabel = {$\alphk$},
			legend style = {legend columns = -1},
		},
		backtrackaxis/.style = {
			myaxis,
			xlabel = {\# of iterations},
			ylabel = {cumulative \# of backtracks},
			legend style = {legend columns = -1},
		},
	}

	\pgfplotsset{
		plotline/.style 2 args = {
			very thick,
			mark = {#1},
			color = {#2},
			mark size = 3pt,
			mark options = {solid},
			mark repeat = 5,
			mark phase = 1,
		},
		plotline*/.style 2 args = {
			plotline = {#1}{#2},
			dashed,
		},
		adabim/.style = {%
			plotline = {o}{black},
		},
		stabim/.style = {%
			plotline* = {x}{black},
		},
		bigsam/.style = {
			plotline = {x}{red},
		},
		solodov/.style = {
			plotline = {o}{green},
		},
		solodov2/.style = {
			plotline = {o}{orange},
		},
		solodov3/.style = {
			plotline = {o}{teal},
		},
		it3D/.style = {
			plotline* = {diamond}{brown},
		},
		bisg/.style = {%
			plotline* = {diamond}{blue},
		},
	}

\makeatother

\newcommand{\TheTitle}{On the convergence of proximal gradient methods\texorpdfstring{\\}{ }for convex simple bilevel optimization%
}
\newcommand{\TheShortTitle}{Proximal gradient methods for convex simple bilevel optimization}
\newcommand{\TheShortAuthors}{P. Latafat, A. Themelis, S. Villa, and P. Patrinos}
\newcommand{\TheKeywords}{%
	Convex optimization%
	\Sep
	bilevel programming%
	\Sep
	adaptive proximal gradient methods%
	\Sep
	locally Lipschitz gradient%
}
\newcommand{\TheSubjclass}{%
	\amsmscLink{65K05}%
	\Sep
	\amsmscLink{90C06}%
	\Sep
	\amsmscLink{90C25}%
	\Sep
	\amsmscLink{90C30}%
}
\newcommand{\TheFunding}{%
	This work was supported by:
	the Research Foundation Flanders (FWO) postdoctoral grant 12Y7622N and research projects G081222N, G033822N, G0A0920N;
	Research Council KU Leuven C1 project No. C14/18/068;
	European Union's Horizon 2020 research and innovation programme under the Marie Skłodowska-Curie grant agreement No. 953348;
	Japan Society for the Promotion of Science (JSPS) KAKENHI grant JP21K17710.
	S. V. acknowledges the support of the European Commission (grant TraDE-OPT 861137), the US Air Force Office of Scientific Research (FA8655-22-1-7034), the Ministry of Education, University and Research (PRIN 202244A7YL project ``Gradient Flows and Non-Smooth Geometric Structures with Applications to Optimization and Machine Learning'').
	The research by S. V. has been supported by the MIUR Excellence Department Project awarded to Dipartimento di Matematica, Università di Genova, CUP D33C23001110001. S. V. is a member of the Gruppo Nazionale per l'Analisi Matematica, la Probabilità e le loro Applicazioni (GNAMPA) of the Istituto Nazionale di Alta Matematica (INdAM).
	This work represents only the view of the authors.
	The European Commission and the other organizations are not responsible for any use that may be made of the information it contains.%
}

\ifarxiv

	\title{\TheTitle\thanks{\TheFunding}}
	\date{}

	\author{%
		Puya Latafat\thanks{\hspace*{0pt}{\TheAddressKU}}%
		\and
		Andreas Themelis\thanks{\hspace*{0pt}{\TheAddressKUJ}}
		\and
		Silvia Villa\thanks{\hspace*{0pt}%
			{\TheAddressUGe.}
			\newline
			{\it E-mails:} \sf
			\{%
				\emailLink[puya.latafat@kuleuven.be]{puya.latafat}%
			,%
				\emailLink[panos.patrinos@kuleuven.be]{panos.patrinos}%
			\}%
			\emailLink[puya.latafat@kuleuven.be,panos.patrinos@kuleuven.be]{@kuleuven.be},
			\emailLink{andreas.themelis@ees.kyushu-u.ac.jp},
			\emailLink{villa@dima.unige.it}%
		}%
		\and
		Panagiotis Patrinos\footnotemark[2]%
	}

	\date{}

\else
	\journalname{JOTA}

	\title{\TheTitle}
	\titlerunning{\TheShortTitle}

	\author{Puya Latafat\and Andreas Themelis\and Silvia Villa\and Panagiotis Patrinos}
	\authorrunning{\TheShortAuthors}

	\institute{%
		Puya Latafat (corresponding author) and Panagiotis Patrinos\at\TheAddressKU\\
		\email{puya.latafat@kuleuven.be, panos.patrinos@esat.kuleuven.be}%
	\and
		Andreas Themelis\at\TheAddressKUJ\\
		\email{andreas.themelis@ees.kyushu-u.ac.jp}%
	\and
		Silvia Villa\at\TheAddressUGe\\
		\email{villa@dima.unige.it}%
	}

	\date{Received: date / Accepted: date}
\fi

\begin{document}

	\maketitle
	\begin{abstract}

		This paper studies proximal gradient iterations for solving simple bilevel optimization problems where both the upper and the lower level cost functions are split as the sum of differentiable and (possibly nonsmooth) proximable functions.
		We develop a novel convergence recipe for iteration varying stepsizes that relies on Barzilai-Borwein type local estimates for the differentiable terms.
		Leveraging the convergence recipe, under global Lipschitz gradient continuity, we establish convergence for a nonadaptive stepsize sequence, without requiring any strong convexity or linesearch.
		In the locally Lipschitz differentiable setting, we develop an adaptive linesearch method that introduces a systematic adaptive scheme enabling large  and nonmonotonic stepsize sequences while being insensitive to the choice of hyperparameters and initialization.
		Numerical simulations are provided showcasing favorable convergence speed of our methods.
	\end{abstract}

	\subclass{\TheSubjclass}

\ifarxiv
\fi

	\section{Introduction}\label{sec:introduction}

		Bilevel programs consist of optimization problems with a hierarchical structure where the solution of one optimization problem is sought over the set of solutions of another one.
		Such problems originally emerged in the framework of game theory and have been studied extensively since the 1950s, see \cite{dempe2002foundations,dempe2020bilevel} for an extensive overview.
		Recently, they have also found applications in various areas of machine learning such as hyperparameter optimization, meta learning, data poisoning attacks, and reinforcement learning \cite{franceschi2018bilevel,rajeswaran2019meta,hong2023two,grazzi2023bilevel,borsos2020coresets}.
		Variational inequality variants have also been of much interest in recent years \cite{facchinei2014vi,bigi2022combining,lampariello2022solution,kaushik2021method,pedregosa2016hyperparameter}.
		The standard approach for addressing bilevel programs consists of solving a series of approximate problems with better regularity properties; refer to \cite{facchinei2014vi,bahraoui1994convergence,attouch1996viscosity} and the references therein.
		However, it is widely known that these techniques can suffer from many practical issues related to convergence speed and stability.

		In this work, we study \emph{simple bilevel programs} which refers to problems where the lower level does not have a parametric dependence on the variables of the upper level problem.
		We split both the upper and the lower cost functions as the sum of differentiable and nonsmooth terms and study two explicit algorithms without the need to solve any inner minimizations.
		In particular, we consider structured simple bilevel programs of the form
		\begin{subequations}\label{eq:P}
			\begin{align}
				\label{eq:P1}
				\minimize_{x\in\R^n}{} &~ \cost_1(x)\coloneqq \f_1(x)+\g_1(x)
			\\
				\label{eq:P2}
				\stt{} &~ x\in\X_2\coloneqq\argmin_{w\in\R^n}\set{\cost_2(w)\coloneqq \f_2(w)+\g_2(w)},
			\end{align}
		\end{subequations}
		where functions \(\f_1,\f_2:\R^n\to \R\) are convex and have (\emph{locally}) Lipschitz continuous gradients, and \(\g_1, \g_2:\R^n\to \Rinf\) are proper closed convex (potentially nonsmooth) functions.
		Some notable example applications include regularized problems in machine learning and signal processing, where the regularization can be captured by the upper level functions, e.g., \(\g_1 = \|\cdot\|_p\), \(p\geq 1\), corresponding to \(\ell^p\) regularization, while the loss function can be captured by the lower level function \(\f_2\), and there may be additional constraints such as nonnegativity constraints captured by \(\g_2\).
		The above formulation also encompasses the class of convex nonlinear programs (NLPs) where nonsmooth proximable terms can be incorporated in the cost function, unlike typical NLP formulations (see \cref{rem:NLP}).

		A fundamental approach for tackling bilevel programs relies on the so-called \emph{diagonal} approach \cite{bahraoui1994convergence,cabot2005proximal,solodov2007explicit}, which involves examining the scaled sum of the upper and lower cost functions
		\begin{equation}\label{eq:phik}
			\fk\coloneqq\sigk\f_1+\f_2,
			\quad\gk\coloneqq\sigk\g_1+\g_2,
		\quad\text{and}\quad
			\phik\coloneqq \fk + \gk = \sigk \cost_1 + \cost_2,
		\end{equation}
		parametrized with a scalar \(\sigk>0\).
		The method by Cabot \cite{cabot2005proximal} for solving simple bilevel programs, here dubbed \emph{Cabot's proximal point algorithm} (\cabot), involves iterative proximal maps (see \cref{sec:preliminaries})
		\begin{equation*}%
			x_{k+1}
		=
			\prox_{\alphk*\phik*}(x_k),
		\end{equation*}
		where the parameter \(\sigk\) is updated after each iteration.
		Convergence of \cabot{} was established under the \emph{slow control} condition
		\begin{equation}\label{eq:sigk}
			\sigk\searrow 0
		\quad\text{and}\quad
			\sum_{k\in\N}\sigk=\infty.
		\end{equation}
		However, due to its implicit nature, in many applications \cabot{} leads to inner minimizations or matrix inversions.
		A notable advancement in this regard was achieved by Solodov in \cite{solodov2007explicit}, who studies \eqref{eq:P} when \(\g_1 \equiv 0\) and \(\g_2\) is an indicator of a closed convex set \(D\).
		The method, here dubbed \emph{Solodov's explicit descent method} (\solodov), uses explicit oracles (gradients for \(\f_1 ,\f_2\) and projections onto \(D\)) and updates \(\sigk\) after a single step of projected gradient method with Armijo linesearch, without the need to solve any inner minimizations.
		More specifically, given \(\nu, \eta \in (0,1)\) and some \(\widehat\alpha_0>0\), in each iteration an inverse penalty \(0<\sigk*\leq\sigk\) is chosen and the variable is updated as
		\begin{subequations}\label{eq:solodov}
			\begin{align}
				x_{k+1}
			={} &
				\proj_D\left(x_k - \alphk*\nabla\fk*(x_k)\right),
			\label{eq:solodov:PG}
			\shortintertext{where \(\alphk*=\widehat\alpha_0\eta^{m_k}\) and \(m_k\in\N\) is the smallest such that}
			\label{eq:solodov:LS}
				\fk*(x_{k+1})
			\leq{} &
				\fk*(x_k)
				+
				\nu
				\innprod{\nabla\fk*(x_k)}{x_{k+1}-x_k}.
			\end{align}
		\end{subequations}
		The challenge, inherent to the bilevel setting, lies in the fact that \(\cost_1\) may take values smaller than the optimal solution along the iterates, making the usual telescoping argument invalid.
		This was overcome in \cite{cabot2005proximal,solodov2007explicit} by an intricate analysis that, while not imposing any strong convexity, does require the set of solutions to be bounded.
		Moreover, the explicit setting of \cite{solodov2007explicit} relies crucially on the fact that \(\g_2\) is an indicator function, making the extension to the proximal case nontrivial.

		In addition to \cite{cabot2005proximal,solodov2007explicit}, several other methods have been proposed for solving simple bilevel problems.
		The \emph{minimal norm gradient} method (\mng) \cite{beck2014first} studies \eqref{eq:P} when \(\g_1\equiv 0\), \(\g_2\) is the indicator of a closed convex set, and the upper level problem is strongly convex.
		\Mng{} relies on a cutting plane approach which can lead to inner minimizations.
		The \emph{bilevel gradient sequential averaging method} (\bigsam) is an explicit method proposed in \cite{sabach2017first} based on a viscosity approximation approach \cite{xu2004viscosity,moudafi2000viscosity}. It considers problems with Lipschitz differentiable and strongly convex upper level cost functions, and establishes an \(O(\nicefrac1K)\) worst-case convergence rate in terms of the lower cost function (see \cref{sec:bigsam} for further details).
		Another related algorithm that we consider in the simulations is the \emph{iterative regularization via dual diagonal descent} (\itthreeD) \cite{garrigos2018iterative}, which is designed for the iterative regularization of linear inverse problems.
		A remarkable property of \itthreeD{} is that it does not impose the slow control condition (see \cite[Rem. 10]{garrigos2018iterative}).
		The \emph{Bi-Sub-Gradient method (version II)} (\bisg) detailed in \cref{sec:bisg} was proposed in \cite{merchav2023convex} that allows for nonsmoothness on both levels.
		It extends \bigsam{} by relaxing the strong convexity assumption. Like \bigsam{}, it achieves a worst-case convergence rate of \(O(\nicefrac1K)\) in terms of the lower level cost function.
		Finally, the \emph{diagonal gradient scheme} (\dgs) was proposed in \cite{peypouquet2012coupling} for solving smooth simple bilevel programs.
		It however involves an implicit stepsize rule that is available in closed form only in certain scenarios, such as when a quadratic growth condition holds (see \cite[Assumptions H1-H3 and \S3.2]{peypouquet2012coupling}).

		The above nonexhaustive literature review largely focused on \emph{explicit} proximal gradient--based methods.
		There have been also many studies that have considered simple bilevel programs with nondifferentiable terms possibly on both levels.
		Most notably, \emph{Solodov's bundle method} (\bundle) proposed in \cite{solodov2007bundle} achieves this through explicit subgradient operations combined with minimization subroutines for constructing cutting-planes approximations.
		Another relevant work is \cite{helou2017subgradient} that offers a unified framework for the analysis of (sub)gradient-type iterations.
		Other recent contributions include \cite{doron2022methodology} that does not require strong convexity or differentiability of the upper level problem, but involves inner subroutines.
		A cutting-planes strategy that employs conditional gradient--type updates is proposed in \cite{jiang2023conditional}.
		In \cite{guan2023first} the authors propose a minimal like--norm gradient method under H\"{o}lderian-type assumptions.

		\subsection{Contributions}%

			We show that proximal gradient iterates involving the family of scaled functions in \eqref{eq:phik} converge under global Lipschitz gradient continuity, without any additional strong convexity assumption.
			Notably, our proposed scheme, \refstabim{}, allows for nonsmoothness in the upper level cost, while only requiring proximal operations for individual nonsmooth terms \(\g_1, \g_2\) (see \cref{rem:prox:gk} for evaluating the proximal mapping).

			The convergence is actually established for a more general framework that identifies three main properties of the stepsize sequence that, combined with the slow control condition \eqref{eq:sigk} are shown to suffice.
			Our analysis crucially relies on utilization of certain Barzilai-Borwein type estimates for the differentiable terms (see \eqref{eq:lL}).
			As a result, diverging from \cite{solodov2007explicit,cabot2005proximal} that establish a quasi-descent inequality in terms of the distance from the solution, we instead show that such a property holds for a more intricate quantity that depends on a combination of cost function, distance from solution(s) and fixed point residual (see \cref{thm:LW}).

			Thanks to this more general perspective, a new linesearch method, \refadabim{}, is proposed, that similar to \solodov{} \cite{solodov2007explicit} can cope with problems involving merely \emph{locally} Lipschitz continuous gradients, while extending the aformentioned work to the proximal setting.
			\refAdabim*{} provides a dynamic update of the initial stepsize \(\balphk\), see \eqref{eq:adabim:alphk*_0}, (as opposed to fixing a predefined stepsize initialization hyperparameter) which is refined over the iterations based on the local Barzilai-Borwein--type estimates, yielding larger stepsizes with considerably fewer backtrackings (see the top row of \cref{fig:gam}).

		\subsection{Organization}%

			After listing some preliminary material, the next section starts by presenting our assumptions along with commentary on their generality, and concludes with the two proposed algorithms.
			The proofs of the methods are deferred to \cref{sec:appendix} and rely on a unifying convergence recipe that is presented in detail in \cref{sec:cnv}.
			After some preliminary lemmas related to the adaptive strategy, we derive a quasi-descent inequality in \cref{sec:prelem};
			this is used in \cref{sec:recipe} for developing the aforementioned convergence recipe.
			Numerical simulations are carried out in \cref{sec:num}, and the paper concludes with some final comments in \cref{sec:conclusions}.

		\subsection{Preliminaries}\label{sec:preliminaries}%

			The sets of natural, real, and extended-real numbers are \(\N\), \(\R\coloneqq(-\infty,\infty)\) and \(\Rinf\coloneqq\R\cup\set\infty\), respectively, while the positive and strictly positive reals are \(\R_+\coloneqq[0,\infty)\) and \(\R_{++}\coloneqq(0,\infty)\).
			We adopt the conventions that \(0\in\N\) and \(1/0=\infty\).
			Given \(x\in\Rinf\), its positive part is indicated as \([x]_+\coloneqq\max\set{0,x}\).
			We denote by \(\innprod{\cdot}{\cdot}\) and \(\|{}\cdot{}\|\) the standard Euclidean inner product and the induced norm, and with \(\id\) the identity function defined on a suitable space.
			The closed Euclidean ball of radius \(r>0\) and centered at \(\bar x\in\R^n\) is denoted as \(\cball{\bar x}{r}\coloneqq\set{x\in\R^n}[\|x-\bar x\|\leq r]\).
			Given two nonempty sets \(U,V\subseteq\R^n\), with \(U+V\coloneqq\set{u+v}[u\in U,v\in V]\) we indicate their Minkowski sum, while \(\conv V\) is used to denote the convex hull of \(V\).


			The \DEF{domain} and \DEF{epigraph} of an extended real--valued function \(\func{h}{\R^n}{\Rinf}\) are, respectively, the sets
			\(
				\dom h
			\coloneqq
				\set{x\in\R^n}[h(x)<\infty]
			\)
			and
			\(
				\epi h
			\coloneqq
				\set{(x, c )\in\R^n\times\R}[h(x)\leq c]
			\).
			Function \(h\) is said to be \DEF{proper} if \(\dom h\neq\emptyset\), and \DEF{lower semicontinuous (lsc)} if \(\epi h\) is a closed subset of \(\R^{n+1}\).
			We say that \(h\) is \emph{level bounded} if its \( c \)-sublevel set
			\(
				\lev_{\leq c }h
			{}\coloneqq{}
				\set{x\in\R^n}[
					h(x)\leq c
				]
			\)
			is bounded for all \( c \in\R\).

			The \DEF{indicator function} of a set $E\subseteq\R^n$ is denoted by \(\indicator_E\), namely \(\indicator_E(x)=0\) if \(x\in E\) and \(\infty\) otherwise.
			The projection onto and the distance from $E$ are respectively denoted by
			\[
				\proj_E(x)\coloneqq\argmin_{z\in E}\|z-x\|
			\quad\text{and}\quad
				\dist(x,E)\coloneqq\inf_{z\in E}\|z-x\|.
			\]
			The (convex) \DEF{subdifferential} of a proper lsc convex function \(\func{h}{\R^n}{\Rinf}\) at a point \(\bar x\) is the set
			\(
				\partial h(\bar x)
			\coloneqq
				\set{v\in\R^n}[h(x)\geq h(\bar x)+\innprod{v}{x-\bar x}\ \forall x\in\R^n]
			\).
			The \DEF{proximal mapping} of \(h\) is \(\func{\prox_h}{\R^n}{\R^n}\) defined by
			\[
				\prox_h(x)
			=
				\argmin_{w\in\R^n}\set{
					h(w)+\tfrac12\|w-x\|^2
				},
			\]
			and is characterized by the implicit subdifferential inclusion \cite[Eq. (24.2)]{bauschke2017convex}
			\begin{equation}\label{eq:proxsubgrad}
				x-\prox_h(x)
			\in
				\partial h\bigl(\prox_h(x)\bigr).
			\end{equation}

	\section{Problem setup and proposed algorithms}\label{sec:ass_results}%

		We will pattern the frameworks of \cite{solodov2007explicit,cabot2005proximal}, while considering  general convex functions \(\g_1\) and \(\g_2\) (as opposed to indicator functions), and differentiable functions \(\f_1, \f_2\) with locally Lipschitz-continuous gradients.
		Our main assumptions are as follows.

		\begin{subequations}\label{eq:costinf}%
			\begin{assumption}[basic requirements]\label{ass:basic}%
				The following hold in problem \eqref{eq:P}:
				\begin{enumeratass}
				\item \label{ass:f}%
					\(\f_1,\f_2:\R^n\to \R\) are convex and have locally Lipschitz-continuous gradients;
				\item \label{ass:g}%
					\(\g_1,\g_2:\R^n \to \Rinf\) are proper lsc convex functions with easy to compute proximal mappings;
				\item \label{ass:cost1}%
					the upper level problem restricted to \(\dom\cost_2\) is lower bounded:
					\begin{equation}\label{eq:costinf_1}
						\costinf_1
					{}\coloneqq{}
						\inf_{x\in\dom\cost_2}\set{\f_1(x)+\g_1(x)}
					{}>{}
						-\infty;
					\end{equation}
				\item \label{ass:X1}%
					the set of solutions \(\X_1\coloneqq\argmin\set{\cost_1(x)}[x\in\argmin\cost_2]\) is nonempty and bounded; in particular,
					\begin{equation}\label{eq:costinf_2}
						\costinf_2
					{}\coloneqq{}
						\inf_{\fillwidthof[c]{x\in\dom \g_2}{x\in\R^n}}\set{\f_2(x)+\g_2(x)}
					{}>{}
						-\infty\mathrlap{.}
					\end{equation}
				\end{enumeratass}
			\end{assumption}
		\end{subequations}

		The boundedness of solution set in \cref{ass:X1} is a standard assumption in the generality of our setting, see \cite{solodov2007explicit}.
		It is in particular implied by conditions such as coercivity of the upper cost function \(\cost_1\). 

		The iterations of our proposed method amount to selecting a stepsize \(\alphk*>0\) together with an (inverse) penalty parameter \(\sigk*\leq\sigk\), followed by one proximal gradient step on the inversely penalized cost function  \(\fk*+\gk*\):
		\begin{equation}\label{eq:PG}
			x_{k+1}
		=
			\prox_{\alphk*\gk*}(x_k-\alphk*\nabla\fk*(x_k)).
		\end{equation}
		\begin{remark}[prox-friendliness of \(\gk\)]\label{rem:prox:gk}
			For the sake of ``explicitness'', we assume throughout that the proximal mapping of \(\gk\) can be evaluated efficiently.
			This can be done without losing generality over \cref{ass:g}, possibly up to lifting of the problem.
			\begin{enumerate}[leftmargin=*,widest=2]
			\item \label{rem:prox:gk:lifted}
				{\em Lifted reformulation.}
				Although in general prox-friendliness is not preserved by the sum, by exploiting the idea presented in \cite{doron2022methodology} the equivalent lifted problem
				\begin{subequations}\label{eq:P:lifted}
					\begin{align}
						\label{eq:P1:lifted}
						\minimize_{x = (z, z')\in\R^{2n}}{} &~ \f_1(z)+\g_1(z)
					\\
						\label{eq:P2:lifted}
						\stt{} &~ x \in \X_2\coloneqq\argmin_{(w,w')\in\R^{2n}}\set{\f_2(w') + \tfrac12\|w-w'\|^2+\g_2(w')}
					\end{align}
				\end{subequations}
				may instead be considered, which by including the quadratic term \(\tfrac12\|w-w'\|^2\) into the smooth component of the lower level results in proximal gradient updates being carried out in parallel (cf. \cref{alg:stabim:gamk*} or \cref{state:adabim:x+}):
				\begin{align*}
					z_{k+1}
				={}&
					\prox_{\alphk*\sigk* \g_1}\left( z_k - \alphk* \bigl(z_k - z_k' + \sigk* \nabla \f_1(z_k)\bigr) \right),
				\\
					z_{k+1}'
				={} &
					\prox_{\alphk* \g_2}\left( z'_k - \alphk* \bigl(z_k' - z_k + \nabla \f_2(z_k')\bigr) \right).
				\end{align*}
				(This reformulation still complies with \cref{ass:basic} provided that \(\f_1+\g_1\) is lower bounded on \(\R^n\), as opposed to merely on \(\dom\g_2\).)
				This procedure can be generalized to address \(\g_1\) of the form \(\g_1(x)=g_1(x)+\dots+g_m(x)\) with each \(g_i\) (individually) being proximable, up to suitably adding slack variables and modifying the smooth term \(\f_2\) in the lower level.

			\item
				{\em Nonlifted option.}
				In many instances of practical interest there is no need to resort to a lifting, as the proximal map of \(\gk\) can be evaluated based on that of \(\g_1\) and/or \(\g_2\).
				This is trivially the case if either \(\g_1\) or \(\g_2\) is zero.
				Other typical instances involve the case in which either one is the \(\ell^1\) norm, or the indicator of a simple set such as a box or an \(\ell^2\)-ball, for which the proximal mapping of the sum is available in closed form; refer to \cite[\S24]{bauschke2017convex} and \cite[\S6]{beck2017first} for examples and further details of proximable functions.
			\qedhere
			\end{enumerate}
		\end{remark}

		\begin{remark}[Nonlinear programs (NLPs)] \label{rem:NLP}
			A notable example is the class of convex NLPs
			\begin{align*}
				\minimize_{x\in\R^n}{} &~ \f_1(x)+\g_1(x)
			\\
				\stt{} &~ x \in D, \quad Ax = b, \quad h(x) \leq 0,
			\end{align*}
			where \(h=(h_1,\dots,h_m)\) and \(\func{\f_1,h_i}{\R^n}{\R}\) are convex and with locally Lipschitz continuous gradients, \(i=1,\dots,m\), \(\g_1:\R^n\to \Rinf\) is a proper closed convex (possibly nonsmooth) proximable function, \(A\in \R^{m\times n}, b\in \R^m\), and \(D\) is a nonempty closed convex set easy to project onto.
			As observed in \cite{solodov2007explicit}, this problem can be formulated in the form of \eqref{eq:P} by setting \(\g_2 = \indicator_D\) and \(\f_2(x) = \|Ax-b\|^2 + \|\max\set{0, h(x)}\|^2\), where \(\indicator_D\) denotes the indicator function of the set \(D\).
			(\Cref{ass:basic} holds provided that the set of solutions is bounded and that \(\inf_D\f_1+\g_1>-\infty\).)
			As commented in \cref{rem:prox:gk:lifted}, as long as \(\f_1+\g_1\) is lower bounded the projection onto \(D\) and the proximal mapping of \(\g_1\) can be decoupled by suitably lifting; further lifting in fact allows for any finite sum structure of proximable terms.
		\end{remark}

		\subsection{The globally Lipschitz case: \Stabim}


			We begin with a linesearch-free proximal gradient method involving the family of scaled functions in \eqref{eq:phik}, under the assumption that \(\f_1\) and \(\f_2\) have globally Lipschitz continuous gradients.
			The full generality of \cref{ass:f} will be addressed in the next subsection through introducing a novel linesearch procedure with an adaptive stepsize initialization.
			The linesearch-free method, \refstabim{}, is \emph{static} in that it uses nonadaptive (and nevertheless increasing) stepsizes (the naming convention will become clear in the sequel).
			It can be viewed as a direct generalization of the iterations considered in \cite{cabot2005proximal,solodov2007explicit} to the full splitting setting of \eqref{eq:P}.
			It extends \cite{cabot2005proximal} by allowing explicit gradient oracles, while it extends \cite{solodov2007explicit} in that it is not limited to constrained problems.

			Although \refstabim*{} uses global estimates and nonadaptive stepsizes, its proof is \emph{only} made possible through the study of a more complex \emph{adaptive stepsize} sequence that allows us to establish a \emph{quasi}-descent behavior on a combination of distance to solutions, cost, and fixed-point residual, see \cref{thm:LW}.
			The main convergence results for \refstabim*{} is presented next.
			Similarly to \cite{solodov2007explicit,cabot2005proximal}, it is shown that the distance from the set of solutions converges to zero.
			Because of nonemptiness and boundedness of the optimal set \(\X_1\) prescribed by \cref{ass:X1}, this condition is equivalent to existence and optimality of the cluster points.

			\begin{theorem}[convergence of \refstabim*]\label{thm:stabim:convergence}%
				Additionally to \cref{ass:basic}, suppose that \(\nabla\f_i\) are globally \(L_{\f_i}\)-Lipschitz continuous, \(i=1,2\), and that \(\seq{\sigk}\) complies with \eqref{eq:sigk}.
				Then (\(\seq{x_k}\) is bounded and) \(\seq{\dist(x_k, \X_1)}\) converges to zero.
				Moreover, the following sublinear rate holds
				\[
					\min_{k\leq K} \|x_{k+1} - x_k\|^2
				\leq
					\tfrac{\nu\bar{\originalvarphi}_0(x_0)}{(1-\nu)L_{\f_2}(1+K)}
				\quad\text{and}\quad
					\min_{k\leq K}\dist^2(0,\partial\phik(x_{k+1}))
				\leq
					\Bigl(1+\sigma_0\tfrac{L_{\f_1}}{L_{\f_2}}\Bigr)
					\tfrac{M_{\rm max}^2\bar{\originalvarphi}_0(x_0)}{(1-\nu)(1+K)},
				\]
				where \(M_{\rm max} = 1+ \nu+\nu\sigma_0\nicefrac{L_{\f_1}}{L_{\f_2}}\) and
				\(
					\bar{\originalvarphi}_0
				\coloneqq
					\sigma_0(\cost_1-\costinf_1) + \cost_2-\costinf_2
				\)
				(\(\costinf_i\) are as in \eqref{eq:costinf}).
			\end{theorem}

			The above worst-case rate result is in line with existing ones in the bilevel setting in terms of the lower cost function, see e.g., \cite{sabach2017first}. Whether convergence rates can be obtained in terms of the upper level cost under additional assumptions is an open problem for future work.

			\begin{algorithm}[t]
			    \caption{%
			        \textbf{\algnamefont Sta}tic \textbf{\algnamefont Bi}level \textbf{\algnamefont M}ethod (\stabim) when \(\protect\f_i\) are globally \(L_{\protect\f_i}\)-Lipschitz smooth, \(i=1,2\)%
			    }%
			    \label{alg:stabim}%

			    	\begin{algorithmic}[1]
    	\itemsep=3pt
    	\Require
    			starting point \(x_0\in\R^n\),~
    			(inverse) penalty \(\sigma_0>0\),~
    			and~
    			\(\nu\in(0,1)\)
    			{\footnotesize (e.g., \(\nu = 0.99\))}

    	\item[\algfont{Repeat for} \(k=0,1,\ldots\) until convergence]

    	\State \label{alg:stabim:sigk*}%
    		Choose \(\sigk*\in\left[\frac{3}{4}\sigk,\sigk\right]\)

    	\State \label{alg:stabim:gamk*}%
    		\(\alphk*=\frac{\nu}{\sigk*L_{\f_1}+L_{\f_2}}\)

    	\State
    		\(
    			x_{k+1}
    		=
    			\prox_{\alphk*\gk*}\left(x_k-\alphk*\nabla\fk*(x_k)\right)
    		\)
    	\Comment{\smash{%
    		\begin{tabular}[t]{@{}l@{}}
    			see \cref{rem:prox:gk} for the proximal evaluation
    		\\
    			of \(\gk*\) based on that of \(\g_1\) and \(\g_2\)
    		\end{tabular}
    	}}%

    	\item[\algfont{Return}]
    		\(x_{k+1}\)
    	\end{algorithmic}
			\end{algorithm}

		\subsection{The locally Lipschitz case: \Adabim}

			The key idea of our adaptive scheme is based on recent works \cite{malitsky2020adaptive,latafat2023adaptive} that study the proximal gradient method, and implicitly enforces a descent inequality bypassing the need for a linesearch.
			While in the bilevel setting a linesearch is still necessary, this analysis provides a systematic adaptive approach for initializing the linesearch.
			Specifically, we use local Lipschitz estimates of the differentiable functions \(\fk=\sigk\f_1+\f_2\) at the previous iterates \(x_k, x_{k-1}\) as%
			\begin{subequations}\label{eq:lL}
				\begin{align}
					\lk
				{}\coloneqq{} &
					\frac{
						\innprod{\nabla\fk(x_{k-1})-\nabla\fk(x_k)}{x_{k-1}-x_k}
					}{
						\|x_{k-1}-x_k\|^2
					}
					,\quad
					\Lk
					{}\coloneqq {}
					\frac{
						\|\nabla\fk(x_{k-1})-\nabla\fk(x_k)\|^2
					}{
						\|x_{k-1}-x_k\|^2
					},
				\shortintertext{%
					and for each \(\f_i\) with
				}
				\label{eq:Lk2}
					\lk^{(i)}
				{}\coloneqq{} &
					\frac{
						\innprod{\nabla\f_i(x_{k-1})-\nabla\f_i(x_k)}{x_{k-1}-x_k}
					}{
						\|x_{k-1}-x_k\|^2
					},
				\quad
					i=1,2,
				\end{align}
				so that \(\lk=\sigk\lk^{(1)}+\lk^{(2)}\).
			\end{subequations}
			Noting that the denominator of \(\lk\), \(\Lk\), or \(\lk^{(i)}\) is zero iff \(\nabla f(x^k)-\nabla f(x^{k-1})=0\), we use the convention \(\nicefrac00=0\) so that both \(\lk\) and \(\Lk\) are (well-defined, positive) real numbers.
			We also remark that these quantities are reminiscent of widely popular Barzilai-Borwein stepsize choices \cite{barzilai1988two} commonly used as heuristics in various minimization settings.
			The proposed \refadabim* uses these quantities in order to establish a quasi-descent inequality (see \cref{thm:LW}) and can be viewed as an extension of \cite[{\algnamefont adaPGM} (Alg. 1)]{latafat2023adaptive};
			if \(\sigk =\sigma\) for all \(k\geq0\) and the linesearch is eliminated, then the algorithm reduces to {\algnamefont adaPGM} applied to the problem of minimizing \(\fk+\gk\).

			The following theorem establishes the convergence results for \refadabim*{} in the full generality of \cref{ass:basic}.
			Its proof relies on the convergence recipe provided in \cref{thm:recipe} and is provided in \cref{sec:appendix}.
			The shown lower bound on the stepsize involves a local Lipschitz modulus for \(\sigma_0\nabla\f_1+\nabla\f_2\) over a compact and convex set \(\V\) that, in addition to containing all the iterates \(x^k\), also includes some of those which were discarded during the linesearch (if any).

			\begin{theorem}[convergence of \refadabim*]\label{thm:adabim:convergence}%
				Suppose that \cref{ass:basic} holds and that \(\seq{\sigk}\) complies with \eqref{eq:sigk}.
				Then, the following holds for the iterates generated by \refadabim{}:
				\begin{enumerate}
				\item \label{thm:adabim:WD}%
					\(
						\alphk
					\geq
						\alpha_{\rm min}
					\coloneqq
						\tfrac{1}{2L_{f_0,\V}}
						\min\set{
							\sqrt{1-\nu}
							,
							\sqrt3\eta\nu
						}
					\)
					for all \(k\geq1\).
					Here, \(L_{f_0,\V}\) is a Lipschitz modulus for \(\nabla f_0=\sigma_0\nabla\f^1+\nabla\f^2\) on the bounded set \(\V\coloneqq\conv\set{x^k}[k\in\N]+\ball{0}{\frac{1-\eta}{\eta}r}\), where
					\(
						r
					\coloneqq
						\max_{k\in\N}\|x^{k+1}-x^k\|
					\).

				\item \label{thm:adabim:sublinear}%
					The following worst-case rates hold
					\[
						\min_{k\leq K} \|x_{k+1} - x_k\|^2
					\leq
						\tfrac{\alpha_{\rm max} \bar{\originalvarphi}_0(x_0)}{(1-\nu)(1+K)}
					\quad\text{and}\quad
						\min_{k\leq K}\dist^2(0,\partial\phik(x_{k+1}))
					\leq
						\tfrac{\alpha_{\rm max}}{\alpha_{\rm min}}
						\tfrac{M_{\rm max}^2\bar{\originalvarphi}_0(x_{0})}{(1-\nu)(1+K)},
					\]
					where \(M_{\rm max} = 1+\alpha_{\rm max} L_{f_0, \V}\) and \(\bar{\originalvarphi}_0\) is as in \cref{thm:stabim:convergence}.

				\item \label{thm:adabim:conv}%
					 (\(\seq{x_k}\) is bounded and)  \(\seq{\dist(x_k, \X_1)}\) converges to zero.
				\end{enumerate}
			\end{theorem}

			We remark that \(\alpha_{\rm max}\) in \cref{state:adabim:x+} is required for theoretical reasons only, while in practice it can be set to a large quantity.
			The only other parameter involved, the initial stepsize \(\alpha_0\), can be set equal to the inverse of a Lipschitz estimate of \(\fk\).
			In referring to \cite[\S2]{latafat2023adaptive} for such practical details,
			we emphasize that the choice of \(\alpha_0\) plays a marginal role, since the update \eqref{eq:adabim:alphk*_0} in combination with the suggested expression of \(\alpha_{-1}\) ensures that any inappropriate initialization is immediately corrected.
			We also remark that our convergence results in \cref{thm:adabim:convergence} are in fact independent of the value of \(\alpha_{-1}\), up to replacing \(\alpha_{\rm min}\gets\min\set{\alpha_{\rm min},\alpha_0}\) in the statement.
			\begin{algorithm}[tb]
					\caption{%
						\textbf{\algnamefont Ada}ptive \textbf{\algnamefont Bi}level \textbf{\algnamefont M}ethod (\adabim) using local Lipschitz estimates \eqref{eq:lL}%
					}%
					\label{alg:adabim}%

			\begin{algorithmic}[1]
			\itemsep=3pt
			\Require
					\begin{tabular}[t]{@{}l@{}}
					starting point \(x_{-1}\in\R^n\),~
					stepsize \(\alpha_0>0\),~
					(inverse) penalty \(\sigma_{-1}=\sigma_0>0\)
				\\
					backtracking linesearch parameters	\(\eta,\nu\in(0,1)\)
					{\footnotesize (e.g., \(\nu = 0.99\), \(\eta = \nicefrac12\))}
				\end{tabular}

			\vspace*{5pt}%
			\Initialize
					\(
						x_0
					{}={}
						\prox_{\alpha_0g_0}\bigl(x_{-1}-\alpha_0\nabla f_0(x_{-1})\bigr)
					\),
					\hfill
					\(\alpha_{\rm max}\gg\frac{1}{\ell_0}\),
					\hfill
					and
					\hfill
					\smash{%
						\(
							\alpha_{-1}
						=
							\alpha_0
							\cdot
							\begin{cases}
								1 & \text{ if }\alpha_0\ell_0\geq\nicefrac{1}{2}\\
								\frac{(\alpha_0\ell_0)^2}{1-(\alpha_0\ell_0)^2} & \text{ otherwise}
							\end{cases}
						\)%
					}%
			\item[\algfont{Repeat for} \(k=0,1,\ldots\) until convergence]
			\State \label{state:adabim:alphk*_0}%
				Choose \(\sigk*\in\left[\frac{3}{4}\sigk,\sigk\right]\),~
				and denoting \(\rhok\coloneqq\frac{\sigk\alphk}{\sigma_{k-1}\alpha_{k-1}}\) initialize the next stepsize as\algfootnotemark\textsuperscript{,}\algfootnotemark
				\begin{equation}\label{eq:adabim:alphk*_0}
					\balphk*
				=
					\tfrac{\sigk}{\sigk*}
					\alphk
					\min\set{
						\sqrt{\tfrac{\sigk}{\sigma_{k-1}}\bigl(1+\rhok\bigr)}
						,\,
						\tfrac{
							\sqrt{
								1
								-
								4\left(1-\frac{\sigk}{\sigma_{k-1}}\right)
								\alphk\lk^{(2)}
							}
						}{
							2\sqrt{[\alphk^2 \Lk^2 - \alphk\lk]_+}
						}
					}
				\end{equation}

			\State \label{state:adabim:x+}%
				\(
					x_{k+1}
				=
					\prox_{\alphk*\gk*}\left(x_k-\alphk*\nabla\fk*(x_k)\right)
				\),
				with \(\alphk*\)
			\Comment{\smash{%
				\begin{tabular}[t]{@{}l@{}}
					see \cref{rem:prox:gk} for the proximal evaluation
				\\
					of \(\gk*\) based on that of \(\g_1\) and \(\g_2\)
				\end{tabular}%
			}}%
			\item[]
				the largest in
				\(\set{ \eta^i\min(\alpha_{\rm max},\balphk*)}[i\in\N]\)
				such that
				\begin{equation}\label{eq:LS}
					\alphk*\lk*\leq\nu
				\end{equation}

			\item[\algfont{Return}]
				\(x_{k+1}\)
			\end{algorithmic}

				\end{algorithm}
				\algfootnotetext{%
					The argument of the square root in the numerator of the second term in \eqref{eq:adabim:alphk*_0} is larger than \(1-\nu >0\) owing to the conditions \(\nicefrac{\sigk*}{\sigk} \geq \nicefrac34\) and \eqref{eq:LS} enforced during the preceding iteration (see \eqref{eq:sqrt>0}).
				}%
				\algfootnotetext{%
					By the convention \(\frac{1}{0}=\infty\), when \(\alphk^2 \Lk^2 - \alphk\lk\leq0\) the definition of \(\balphk*\) reduces to the first term in the minimum.%
				}%

			\subsubsection{Observations about the stepsizes}
				We make some observations about the stepsize sequence of \refadabim*.
				In the initialization at \cref{state:adabim:alphk*_0}, when \(\alphk^2\Lk^2 - \alphk \lk\leq 0\), \ie,
				\begin{equation}\label{eq:ck}
					\ck
				\coloneqq
					\frac{\Lk^2}{\lk}
				=
					\frac{
						\|\nabla\fk(x_{k-1})-\nabla\fk(x_k)\|^2
					}{
						\innprod{\nabla\fk(x_{k-1})-\nabla\fk(x_k)}{x_{k-1}-x_k}
					}
				\leq
					\frac{1}{\alphk},
				\end{equation}
				the second term in \eqref{eq:adabim:alphk*_0} reduces to \(\nicefrac10=\infty\) and the stepsize initialization simplifies as
				\(
					\balphk*
				=
					\sqrt{\tfrac{\sigk}{\sigma_{k-1}}\bigl(1+\rhok\bigr)}
					\tfrac{\sigk}{\sigk*}
					\alphk
				\).
				This is a critical feature since it allows \(\balphk*\) to strictly increase compared to the stepsize \(\alphk\); for instance, under the standard choice \(\sigk = \nicefrac1{k+1}\), the first term is always larger than \(\alphk\).
				There is an apparent trade-off between the first and the second term: a large first term allows for faster recovery from small stepsizes at the expense of a smaller second one which affects the global lower bound on \(\alphk\).
				While this interplay can be tweaked by introducing additional algorithmic parameters as done in \cite{latafat2023convergence}, for clarity of exposition we limit the discussion to this simpler setting.

				\begin{figure}[t]
					\centering
					\includetikz[width=.9\linewidth]{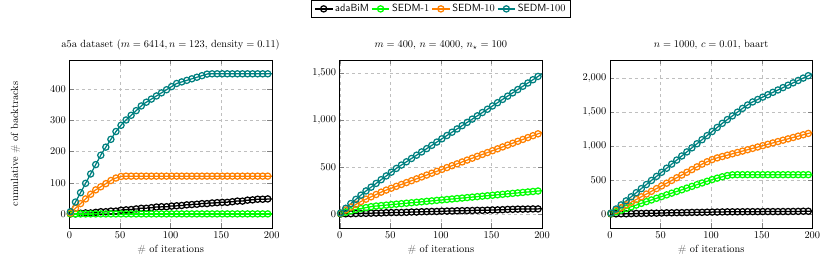}%
					\includetikz[width=.9\linewidth]{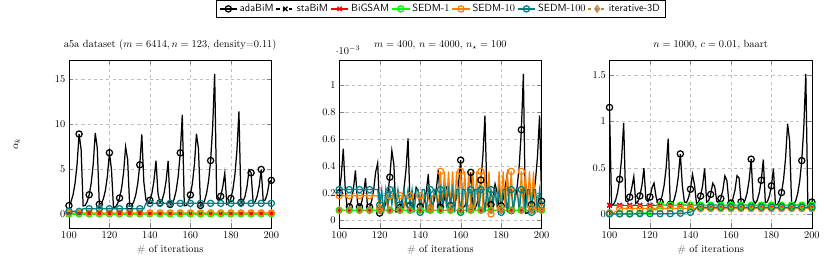}%
					\caption[]{%
						A representative plot showing
						cumulative number of backtracks needed by \refadabim* and \solodov{} (top row)
						and
						stepsize magnitudes in a window of 100 iterations (bottom row)
						in sample simulations from \cref{sec:num}.
						The numerical suffixes {\algnamefont-1}, {\algnamefont-10}, and {\algnamefont-100} in \solodov{} indicate different choices for the value of \(\widehat\alpha_0\) as defined in \cref{sec:solodov}.
						Left: logistic regression (a5a dataset);
						center: linear inverse problem;
						right: solution of integral equations.%
					}%
					\label{fig:gam}%
				\end{figure}

				\refAdabim*{} also retains the linesearch nature of \solodov{}, but compared to the fixed value \(\widehat\alpha_0\) of \eqref{eq:solodov} it provides a dynamic update of the initial stepsize \(\balphk\), cf. \eqref{eq:adabim:alphk*_0}, which is refined over the iterations based on local estimates in \eqref{eq:lL}, yielding much larger stepsizes with considerably fewer backtrackings.
				This online self-correcting feature renders the proposed algorithm insensitive to parameters chosen at initalization.
				As the numerical evidence in \cref{sec:num} well highlights, the overhead caused by the backtracking is negligible compared to the speedup that the dynamic update of \refadabim*{} yields, and even under global Lipschitzian assumption this method exhibits superior performance compared to the \emph{static} counterpart \refstabim*{} presented in the previous subsection.

				The top row of \cref{fig:gam} illustrates the cumulative number of backtrackings per iteration required by \refadabim*{} compared to \solodov{} with different choices of \(\widehat\alpha_0\); high values of \(\widehat\alpha_0\) allow for large stepsizes to be tested and potentially accepted, thereby favoring convergence speed in terms of number of iterations, but may lead to more backtrackings and function evaluations in the linesearch \eqref{eq:solodov:LS}.
				Conversely, small values of \(\widehat\alpha_0\) reduce the complexity of each iteration by reducing the number of backtrackings at the expense of smaller stepsizes and consequently slower convergence.
				These plots correspond to the respective simulations of \cref{sec:numsim} where comparisons in terms of the total number of gradient evaluations are also presented.
				We also remark that the oscillatory behavior of the stepsizes of \refadabim* reported in the bottom plots of \cref{fig:gam} is a key element enabling its fast convergence and has been observed also in the setting of minimization problems, cf. \cite[\S4.3]{latafat2023adaptive}.

	\section{Convergence analysis}\label{sec:cnv}

			In this section we examine the convergence properties of \refadabim* (resp. \refstabim*) for solving problem \eqref{eq:P} under local (resp. global) Lipschitz continuity of the gradients of \(\f_1\) and \(\f_2\).
			Crucially, regardless of the stepsize selection strategy, our analysis relies on a quasi-descent inequality for the proximal gradient updates \eqref{eq:PG}.
			This allows for a unified analysis provided in \cref{thm:recipe} that is based on the identification of a set of properties that the stepsizes should satisfy in order to guarantee convergence.
			We begin by elaborating on some of the notational conventions; a full list is synopsized in \cref{table:notation} for convenience.

			\begin{remark}[Bar notation for minima and shifted costs]\label{rem:cost}%
				Let us define
				\[
					\bcost_i\coloneqq\cost_i-\costinf_i, \; i=1,2,
				\quad\text{and}\quad
					\bphik \coloneqq \sigk \bcost_1 + \bcost_2,
				\]
				where \(\costinf_i=\inf_{\dom\cost_2}\cost_i\) as in \cref{ass:basic} (restricting to \(\dom\cost_2\) is superfluous for \(i=2\)), and let
				\begin{equation}\label{eq:phi*}
					\costinf*
				\coloneqq
					\inf_{\X_2}\cost_1
				\end{equation}
				be the optimal cost of problem \eqref{eq:P}.
				Then,
				\begin{enumerate}
				\item\label{thm:bar>=0}%
					\(\bcost_1(x)\geq0\) for any \(x\in\dom\cost_2\) and \(\bcost_2(x)\geq0\) for any \(x\in\R^n\);
				\item \label{thm:Phi*}%
					\( \tfrac1{\sigk}\left(\phik(x_\star)-\costinf_2\right)=\costinf*\) for any \(x_\star\in\X_1\) and \(k\in\N\).
				\qedhere
				\end{enumerate}
			\end{remark}

			A key step in our convergence analysis based on adaptive stepsizes relies on  the introduction of the quantities in \eqref{eq:lL}.
			We define the shorhand notation for the forward operator \(\Hk = \id - \alphk\nabla\fk\), and note that by optimality conditions of the prox-grad update \eqref{eq:PG}
			\begin{equation}\label{eq:Hk:phi}
				\tfrac1{\alphk}(\Hk(x_{k-1})-\Hk(x_k)) \in \partial \phik(x_k).
			\end{equation}
			As we state in the next lemma, not only do the quantities in \eqref{eq:lL} provide an exact description of the local Lipschitz modulus of \(\nabla \fk\), but also that of the forward operator \(\Hk\).

			\begin{fact}[{\cite[Lem. 2.1]{latafat2023adaptive}}]\label{thm:lL}%
				Suppose that \cref{ass:f} holds, and for \(x_{k-1},x_k\in\R^n\) and \(\alphk>0\) let \(\lk\) and \(\Lk\) be as in \eqref{eq:lL} and \(\Hk\coloneqq\id-\alphk\nabla\fk\).
				Then, the following hold:
				\begin{enumerate}
				\item \label{thm:H}%
					\(
						\|\Hk(x_{k-1})-\Hk(x_k)\|
					=
						\Mk
						\|x_{k-1}-x_k\|
					\),
					where
					\[
						|1-\alphk\Lk|
					\leq
						\Mk
					\coloneqq
						\sqrt{1 + \alphk^2 \Lk^2 -2\alphk\lk}
					\leq
						1+\alphk\Lk.
					\]
				\item\label{thm:l<=L}%
					\(
						\lk
					\leq
						\Lk
					\leq
						\sigk L_{\f_1,\V} + L_{\f_2,\V}
					\leq
						\sigma_0 L_{\f_1,\V} + L_{\f_2,\V}
					\),
					where \(L_{\f_i,\V}\) is a Lipschitz modulus for \(\nabla\f_i\) on a convex set \(\V\) containing \(x_{k-1}\) and \(x_k\), \(i=1,2\).\footnote{%
						The first inequality owes to Cauchy-Schwarz; the last one to convexity of \(\f_1\) together with the fact that \(\sigk\leq\sigma_0\).%
					}%
				\end{enumerate}
			\end{fact}


		\phantomsection\addcontentsline{toc}{subsection}{Notational conventions}\label{sec:notation}%
		\subsection*{Notational conventions}

			As done above and throughout, we use subscripts for iteration counters, typically \(k\), and bracketed superscripts to indicate the level (either 1 or 2).
			Other symbols that will be introduced for the sake of the convergence proofs adhere to the same conventions and are collected in \cref{table:notation}, inclusive of references to the respective definitions (those which are local to the scope of individual proofs are omitted from the list).
			In particular, the uppercase \(\Fk\), \(\Gk\) and \(\Phik\) will be useful for the convergence analysis, and correspond to the respective lowercase symbols scaled by \(\frac{1}{\sigk}\).
			Keeping in mind that \(\sigk\) is an \emph{inverse} penalty parameter, in the sense that it is driven to 0, we refer to \(\Phik=\cost_1+\frac{1}{\sigk}\cost_2\) as the penalized cost, and to \(\phik=\sigk\cost_1+\cost_2\) as the \emph{inversely} penalized cost of the single-level subproblems.

			\begin{table}[htb]
				\centering
				\footnotesize
				\let\mr\multirow
				\renewcommand{\arraystretch}{1.35}%
				\setlength{\tabcolsep}{2.5pt}%
				\begin{tabular}[t]{@{}|c@{~~}l|c| c |c@{~~}l|c|@{}}%
					\multicolumn{3}{@{}c}{\bf upper level \(\cost_1=\f_1+\g_1\)}
					& \multicolumn{1}{c}{} &
					\multicolumn{3}{c@{}}{\bf lower level \(\cost_2=\f_2+\g_2\)}
				\\\cline{1-3}\cline{5-7}
						\(\f_1\)  & smooth part     & \mr{2}{*}{\eqref{eq:P1}}
						&&
						\(\f_2\)  & smooth part     & \mr{2}{*}{\eqref{eq:P2}}
				\\
						\(\g_1\)  & proximable part &
						&&
						\(\g_2\)  & proximable part &
				\\\cline{1-3}\cline{5-7}
						\(\costinf_1\) & \(\inf_{\dom\cost_2}\cost_1\) & \eqref{eq:costinf_1}
						&&
						\(\costinf_2\) & \(\inf\cost_2\)            & \eqref{eq:costinf_2}
				\\\cline{1-3}\cline{5-7}
						\(\bcost_1\) & \(\cost_1-\costinf_1\)  (\(\geq0\) on \(\dom\cost_2\)) & \cref{rem:cost}
						&&
						\(\bcost_2\) & \(\cost_2-\costinf_2\)  (\(\geq0\))                    & \cref{rem:cost}
				\\\cline{1-3}\cline{5-7}
						\(\X_1\) & \(\argmin_{\X_2}\cost_1\)   (optimal set)                  & \cref{ass:X1}
						&&
						\(\X_2\) & \(\argmin\cost_2\)          (feasible set)                 & \eqref{eq:P2}
				\\\cline{1-3}\cline{5-7}
						\(\costinf*\) & \(\min_{\X_2}\cost_1\) (optimal cost)                 & \eqref{eq:phi*}
				\\\cline{1-3}
				\multicolumn{7}{c}{}\\[-2.5ex]
					\multicolumn{3}{@{}c}{\bf single-level inverse-penalty reformulation}
					& \multicolumn{1}{c}{} &
					\multicolumn{3}{@{}c}{\bf single-level penalty reformulation}
				\\\cline{1-3}\cline{5-7}
						\(\fk\)  & \(\sigk\f_1+\f_2\) smooth part               & \mr{3}{*}{\eqref{eq:phik}}
						&&
						\(\Fk\)  & \(\f_1+\frac{1}{\sigk}\f_2\) smooth part     & \mr{3}{*}{\eqref{eq:uppercase}}
				\\
						\(\gk\)  & \(\sigk\g_1+\g_2\) proximable part           &
						&&
						\(\Gk\)  & \(\g_1+\frac{1}{\sigk}\g_2\) proximable part &
				\\
						\(\phik\) & \(\fk+\gk=\sigk\cost_1+\cost_2\)            &
						&&
						\(\Phik\)  & \(\Fk+\Gk=\cost_1+\frac{1}{\sigk}\cost_2\) &
				\\\cline{1-3}\cline{5-7}
						\(\bphik\) & \(\sigk\bcost_1+\bcost_2\) (\(\geq0\) on \(\dom\cost_2\))                   & \cref{rem:cost}
						&&
						\(\Delk\)        & \(\nicefrac{1}{\sigk}-\nicefrac{1}{\sigma_{k-1}}\)                    & \cref{thm:LW}
				\\\cline{1-3}\cline{5-7}
				\multicolumn{7}{c}{}\\[-2.5ex]
					\multicolumn{3}{@{}c}{\bf algorithmic parameters}
					& \multicolumn{1}{c}{} &
					\multicolumn{3}{@{}c}{\bf adaptive estimates}
				\\\cline{1-3}\cline{5-7}
						\(\alphk*\)   & stepsize                                          & \mr{3}{*}{\cref{alg:adabim}}
						&&
						\(\lk\)       & \mr{2}{*}{Lipschitz estimates of \(\nabla\fk\) at \(x_k\)} & \mr{3}{*}{\eqref{eq:lL}}
				\\\cline{1-2}
						\(\sigk*\)    & (inverse) penalty                                 &
						&&
						\(\Lk\)       &                                                   &
				\\\cline{1-2}\cline{5-6}
						\(\rhok*\)    & \(\nicefrac{\sigk*\alphk*}{\sigk\alphk}\)         &
						&&
						\(\lk^{(i)}\) & a Lipschitz estimate of \(\nabla\f_i\) at \(x_k\)   &
				\\\cline{1-3}\cline{5-7}
				\end{tabular}
				\caption{%
					Schematics of the notation adopted in the paper with references to their definitions.
				}%
				\label{table:notation}%
			\end{table}


		\subsection{A quasi-descent inequality}\label{sec:prelem}

			Before delving into the convergence analysis, we will present a series of preliminary results that can be regarded as an extension of the adaptive mechanism of {\algnamefont adaPGM} proposed in \cite[Alg. 1]{latafat2023adaptive}.
			This adaptive scheme not only significantly improves the computational efficiency of our approach, but it also allows us to consider nonsmooth terms on both levels. 

			We proceed to investigate the progress of a single proximal gradient step as described in \eqref{eq:PG} using an arbitrary stepsize \(\alphk*>0\).
			Departing from \cite{solodov2007explicit,cabot2005proximal}, the key of our convergence analysis, captured in the following lemma, is the adoption of penalized (as opposed to inversely penalized) costs.
			As already mentioned in the preview of \cref{table:notation}, we adopt an uppercase notation for the penalized cost
			\begin{equation}\label{eq:uppercase}
				\Phik=\Fk+\Gk
			\quad\text{with}\quad
				\Fk=\tfrac1{\sigk}\fk=\f_1+\tfrac1{\sigk}\f_2
			\quad\text{and}\quad
				\Gk=\tfrac1{\sigk}\gk=\g_1+\tfrac1{\sigk}\g_2,
			\end{equation}
			noticing that \(\Phik=\cost_1+\frac{1}{\sigk}\cost_2=\frac{1}{\sigk}\phik\).
			Doing so allows us to express the difference
			\begin{equation}
				\Phik*-\Phik
			{}={}
				\left(\tfrac1{\sigk*}-\tfrac1{\sigk}\right)\cost_2
			\end{equation}
			as a multiple of the lower-level cost \(\cost_2\), rather than of the upper-level cost \(\cost_1\).

			\begin{lemma}[quasi-descent inequality]\label{thm:LW}%
				Suppose that \cref{ass:basic} holds and consider iterations \eqref{eq:PG} with \(0<\sigk*\leq\sigk\).
				Let \(\bcost_2\coloneqq\cost_2-\costinf_2\geq0\) be as in \cref{rem:cost}, and define \(\Delk\coloneqq\frac{1}{\sigk}-\frac{1}{\sigma_{k-1}}\), \(\rhok*\coloneqq\frac{\sigk*\alphk*}{\sigk\alphk}\), and
				\begin{align*}
					\mathcal L_k(x_\star)
				\coloneqq{} &
					\tfrac12\|x_k-x_\star\|^2
					+
					W_k
				\shortintertext{with}
					W_k
				\coloneqq{} &
					\tfrac{1-4\alphk\sigk\lk^{(2)}\Delk}{4}\|x_k-x_{k-1}\|^2
					+
					\tfrac{\sigk*}{\sigk}\alphk*\rhok*
					\bcost_2(x_{k-1})
				\\
				&
					+
					\sigk\alphk\Delk\bcost_2(x_k)
					+
					\sigk\alphk\left(1+\rhok\right)
					\bigl(
						\cost_1(x_{k-1})
						-
						\costinf*
					\bigr).
				\end{align*}
				Then, for every \(k\in\N\) and \(x_\star\in\X_1\) it holds that
				\begin{align*}
					\mathcal L_{k+1}(x_\star)
					-
					\mathcal L_k(x_\star)
				\leq{} &
					-
					\left(
						\tfrac{1}{4}
						+
						\rhok*^2(\alphk\lk-\alphk^2\Lk^2)
						-
						\sigk\alphk\lk^{(2)}\Delk
					\right)
					\|x_{k-1}-x_k\|^2
				\\
				&
					-
					\sigk\alphk
					\bigl(
						1+\rhok-\rhok*^2
					\bigr)
					\bigl(\cost_1(x_{k-1})-\costinf*\bigr)
				\\
				&
					-
					\sigk*\alphk*
					\left(
						\tfrac{\Delk}{\rhok*}
						+
						\tfrac{1}{\sigk}\left(
							1+\rhok*
							-
							\rho_{k+2}^2\tfrac{\sigk}{\sigk*}
						\right)
					\right)
					\bcost_2(x_k).
				\numberthis\label{eq:puya:muk:descent}
				\end{align*}
				\begin{proof}
				Let
				\begin{equation}
					p_k%
				\coloneqq
					\phik(x_k)-\phik(x_\star)
				=
					\sigk\bigl(\cost_1(x_k)-\costinf*\bigr)
					+
					\bcost_2(x_k)
				\end{equation}
				(which is not necessarily positive).
				We will prove the claim using the uppercase notation of \eqref{eq:uppercase} and, consistently, set
					\(P_k%
					\coloneqq{\tfrac{1}{\sigk}}p_k%
					=\Phik(x_k)-\Phik(x_\star)
					\).
					Being a simple matter of multiplicative constants, note that proximal gradient iterations \eqref{eq:PG} can equivalently be expressed in the penalized cost reformulation \(\Phik*=\Fk*+\Gk*\) up to suitable scaling of the stepsize, namely,
					\[
						x_{k+1}
					=
						\prox_{\sigk*\alphk*\Gk*}(x_k-\sigk*\alphk*\nabla\Fk*(x_k)).
					\]
					The subgradient characterization of the proximal mapping as in \eqref{eq:Hk:phi} yields
					\begin{equation}\label{eq:subgrad}
						\tfrac{\Hk(x_{k-1})-x_k}{\sigk\alphk}
					=
						\tfrac{x_{k-1}-x_k}{\sigk\alphk}-\nabla\Fk(x_{k-1})
					\in
						\partial\Gk(x_k),
					\end{equation}
					where we remind that \(\Hk=\id-\alphk\nabla\fk\), cf. \eqref{eq:Hk:phi}.
					Hence, since \(\partial\Phik=\nabla\Fk+\partial\Gk\),
					\begin{subequations}
						\begin{align}
						\nonumber
							0
						\leq{} &
							\Phik(x_{k-1})
							-
							\Phik(x_k)
							-
							\tfrac{1}{\sigk\alphk}
							\innprod{x_{k-1}-x_k}{x_{k-1}-x_k}
							+
							\innprod{\nabla\Fk(x_{k-1})-\nabla\Fk(x_k)}{x_{k-1}-x_k}
						\\
						\nonumber
						={} &
							\Phik(x_{k-1})
							-
							\Phik(x_k)
							-
							\tfrac{1}{\sigk\alphk}
							\|x_{k-1}-x_k\|^2
							+
							\tfrac1{\sigk}\innprod{\nabla\fk(x_{k-1})-\nabla\fk(x_k)}{x_{k-1}-x_k}
						\\
						={} &
							\Phik(x_{k-1})
							-
							\Phik(x_k)
							-
							\tfrac{1-\alphk\lk}{\sigk\alphk}
							\|x_{k-1}-x_k\|^2
						\label{eq:phikdescent}
						\\
						={} &
							P_{k-1}%
							-
							P_k%
							+
							(\tfrac{1}{\sigk}-\tfrac{1}{\sigma_{k-1}})\bcost_2(x_{k-1})
							-
							\tfrac{1-\alphk\lk}{\sigk\alphk}\|x_{k-1}-x_k\|^2.
						\label{eq:ineqsubgrad}
						\end{align}
					\end{subequations}
					Again from \eqref{eq:subgrad}, this time with \(k\gets k+1\), we have
					\begin{subequations}\label{subeq:3ineq}
						{\mathtight[0.6]%
						\begin{align*}
							0
						\leq{} &
							\Gk*(x_\star)
							-
							\Gk*(x_{k+1})
							+
							\innprod{\nabla\Fk*(x_k)}{x_\star-x_{k+1}}
							-
							\tfrac{1}{\sigk*\alphk*}
							\innprod{x_k-x_{k+1}}{x_\star-x_{k+1}}
						\\
						={} &
							\Gk*(x_\star)
							-
							\Gk*(x_{k+1})
							+
							\innprod{\nabla\Fk*(x_k)}{x_\star-x_{k+1}}
							+
							\tfrac{1}{2\sigk*\alphk*}
							\|x_k-x_\star\|^2
							-
							\tfrac{1}{2\sigk*\alphk*}
							\|x_\star-x_{k+1}\|^2
						\\
						&
							-
							\tfrac{1}{2\sigk*\alphk*}
							\|x_k-x_{k+1}\|^2
						\\
						={} &
							\Gk*(x_\star)
							-
							\Gk*(x_{k+1})
							+
							\innprod{\nabla\Fk*(x_k)}{x_\star-x_k}
							+
							\innprod{\nabla\Fk*(x_k)}{x_k-x_{k+1}}
							+
							\tfrac{1}{2\sigk*\alphk*}
							\|x_k-x_\star\|^2
						\\
						&
							-
							\tfrac{1}{2\sigk*\alphk*}
							\|x_\star-x_{k+1}\|^2
							-
							\tfrac{1}{2\sigk*\alphk*}
							\|x_k-x_{k+1}\|^2
						\\
						\leq{} &
							\Gk*(x_\star)
							-
							\Gk*(x_{k+1})
							+
							\fillwidthof[c]{
								\innprod{\nabla\Fk*(x_k)}{x_\star-x_k}
							}{
								\Fk*(x_\star)-\Fk*(x_k)
							}
							+
							\underbracket*[0.5pt]{
								\innprod{\nabla\Fk*(x_k)}{x_k-x_{k+1}}
							}_{\text{(A)}}
							+
							\tfrac{1}{2\sigk*\alphk*}
							\|x_k-x_\star\|^2
						\\
						&
							-
							\tfrac{1}{2\sigk*\alphk*}
							\|x_\star-x_{k+1}\|^2
							-
							\tfrac{1}{2\sigk*\alphk*}
							\|x_k-x_{k+1}\|^2,
						\numberthis\label{eq:ineq1}
						\end{align*}}%
						where the last inequality uses convexity of \(\Fk*\).
						As to term (A), we have
						\begin{align*}
							\text{(A)}
						={} &
							\tfrac{1}{\sigk\alphk}
							\innprod{\Hk(x_{k-1})-x_k}{x_{k+1}-x_k}
							+
							\tfrac{1}{\sigk\alphk}
							\innprod{\Hk(x_{k-1})-x_k+\sigk\alphk\nabla\Fk*(x_k)}{x_k-x_{k+1}}
						\\
						\overrel[\leq]{\eqref{eq:subgrad}}{} &
							\fillwidthof[c]{
								\tfrac{1}{\sigk\alphk}
								\innprod{\Hk(x_{k-1})-x_k}{x_{k+1}-x_k}
							}{
								\Gk(x_{k+1})-\Gk(x_k)
							}
							+
							\innprod{\nabla\Fk*(x_k)-\nabla\Fk(x_k)}{x_k-x_{k+1}}
						\\
						&
							+
							\underbracket[0.5pt]{
								\tfrac{1}{\sigk\alphk}
								\innprod{\Hk(x_{k-1})-\Hk(x_k)}{x_k-x_{k+1}}
							}_{\text{(B)}}.
						\numberthis\label{eq:ineqA}
						\end{align*}
						Next, we bound the term (B) by \(\epsk*\)-Young's inequality as
						\begin{align*}
							\text{(B)}
						\leq{} &
							\tfrac{\epsk*}{2\sigk\alphk}\|x_k-x_{k+1}\|^2
							+
							\tfrac{1}{2\epsk*\sigk\alphk}\|\Hk(x_{k-1})-\Hk(x_k)\|^2
						\\
						\overrel{\ref{thm:H}}{} &
							\tfrac{\epsk*}{2\sigk\alphk}\|x_k-x_{k+1}\|^2
							+
							\tfrac{\Mk^2}{2\epsk*\sigk\alphk}\|x_{k-1}-x_k\|^2.
						\numberthis\label{eq:ineqB}
						\end{align*}
					\end{subequations}
					Combining the three inequalities \eqref{subeq:3ineq} yields
					{\mathtight[0.9]%
					\begin{align*}
						0
					\leq{} &
						\Phik*(x_\star)
						-
						\Gk*(x_{k+1})
						-
						\Fk*(x_k)
						+
						\Gk(x_{k+1})
						-
						\Gk(x_k)
						+
						(\tfrac1{\sigk*}-\tfrac{1}{\sigk})
						\innprod{\nabla \f_2(x_k)}{x_k-x_{k+1}}
					\\
					&
						+
						\Biggl\{
							\tfrac{\Mk^2}{2\epsk*\sigk\alphk}\|x_{k-1}-x_k\|^2
							-
							\left(
								\tfrac{1}{2\sigk*\alphk*}
								-
								\tfrac{\epsk*}{2\sigk\alphk}
							\right)
							\|x_k-x_{k+1}\|^2
							-
							\tfrac{1}{2\sigk*\alphk*}
							\|x_\star-x_{k+1}\|^2
					\\
					&
						\hphantom{+\Biggl\{{}}
							+
							\tfrac{1}{2\sigk*\alphk*}
							\|x_k-x_\star\|^2
						\Biggr\}
					\\
					\overrel{\ref{thm:Phi*}}{} &
						\costinf*
						+
						\tfrac1{\sigk*}\costinf_2
						-
						\Fk(x_k)
						-
						\Gk(x_k)
						+
						\bigl(\tfrac1{\sigk*}-\tfrac{1}{\sigk}\bigr)
						\vphantom{\f^2}
						\left(
							\innprod{\nabla \f_2(x_k)}{x_k-x_{k+1}}
							-
							\f_2(x_k)
							-
							\g_2(x_{k+1})
						\right)
					\\
					&
						+
						\{\cdots\}
					\\
					={} &
						-P_k%
						+
						\overbracket[0.5pt]{
							\left(
								\tfrac1{\sigk*}-\tfrac{1}{\sigk}
								\vphantom{\f^2}
							\right)
						}^{\geq0}
						\,
						\overbracket[0.5pt]{
							\left(
								\innprod{\nabla \f_2(x_k)}{x_k-x_{k+1}}
								-
								\f_2(x_k)
								-
								\g_2(x_{k+1})
								+
								\costinf_2
							\right)
						}^{\text{(D)}}
						+
						\Bigl\{\cdots\Bigr\}.
					\end{align*}}%
					By using convexity of \(\f_2\) we can bound the term (D) as
					{\mathtight[0.75]%
						\begin{align*}
							\text{(D)}
						={} &
							\innprod{\nabla \f_2(x_{k+1})}{x_k-x_{k+1}}
							-
							\f_2(x_k)
							-
							\g_2(x_{k+1})
							+
							\innprod{\nabla \f_2(x_{k+1})-\nabla \f_2(x_k)}{x_{k+1}-x_k}
							+
							\costinf_2
						\\
						\leq{} &
							\fillwidthof[c]{
								\innprod{\nabla \f_2(x_{k+1})}{x_k-x_{k+1}}
								-
								\f_2(x_k)
							}{
								-\f_2(x_{k+1})
							}
							-
							\g_2(x_{k+1})
							+
							\fillwidthof[c]{
								\innprod{\nabla \f_2(x_{k+1})-\nabla \f_2(x_k)}{x_{k+1}-x_k}
							}{
								\lk*^{(2)}\|x_{k+1}-x_k\|^2
							}
							+
							\costinf_2,
						\end{align*}
					}%
					which plugged in the previous inequality results in
					{\mathtight[0.75]%
					\begin{align*}
						0
					\leq{} &
						-P_k%
						-
						(\tfrac1{\sigk*}-\tfrac{1}{\sigk})\bcost_2(x_{k+1})
						+
						(\tfrac1{\sigk*}-\tfrac{1}{\sigk})
						\lk*^{(2)}\|x_{k+1}-x_k\|^2
						+
						\tfrac{\Mk^2}{2\epsk*\sigk\alphk}\|x_{k-1}-x_k\|^2
					\\
					&
						-
						\left(
							\tfrac{1}{2\sigk*\alphk*}
							-
							\tfrac{\epsk*}{2\sigk\alphk}
						\right)
						\|x_k-x_{k+1}\|^2
						-
						\tfrac{1}{2\sigk*\alphk*}
						\|x_\star-x_{k+1}\|^2
						+
						\tfrac{1}{2\sigk*\alphk*}
						\|x_k-x_\star\|^2.
					\numberthis\label{eq:3ineqs}
					\end{align*}}%
					Summing \eqref{eq:3ineqs}\(+\betk*\)\eqref{eq:ineqsubgrad}, multiplying by \(\sigk*\alphk*\), and rearranging yields that for every \(\betk*\geq0\) and \(\epsk*>0\) the following hold.
				\begin{align*}
					& \tfrac{1}{2}
					\|x_{k+1}-x_\star\|^2
					+
					\tfrac{\sigk*}{\sigk}\alphk*(1+\betk*)p_k%
					+
					\alphk*\sigk*\Delk*\bcost_2(x_{k+1})
				\\
					+ {}&
					\left(
						\tfrac{1-\rhok*\epsk*}{2}
						-
						\alphk*\sigk*\Delk*%
						\lk*^{(2)}
					\right)
					\|x_k-x_{k+1}\|^2
				\leq
					\tfrac{1}{2}
					\|x_k-x_\star\|^2
					+
					\tfrac{\sigk*}{\sigma_{k-1}}\betk*\alphk*p_{k-1}%
				\\
					-{}&
					\rhok*
					\left(
						\betk*(1-\alphk\lk)
						-
						\tfrac{\Mk^2}{2\epsk*}
					\right)
					\|x_{k-1}-x_k\|^2
					+
					\alphk*\sigk*\betk*\Delk\bcost_2(x_{k-1}).
				\numberthis
				\end{align*}
				By selecting \(\betk=\rhok\) and \(\epsk\coloneqq\nicefrac{1}{2\rhok}\) the inequality reduces to the claimed quasi-descent in terms of \(\mathcal L_{k}(x^\star)\).
				\end{proof}
			\end{lemma}

			By looking at the update rule \eqref{eq:adabim:alphk*_0} for \(\balphk*\), it is apparent that the choice of stepsizes in \cref{alg:adabim} is designed so as to ensure that all the multiplying coefficients on the right-hand side of \eqref{eq:puya:muk:descent} are positive.
			Even so, it should be noted that the inequality does not, in general, imply a monotonic decrease of \(\mathcal L_k(x_\star)\) along the iterates, the reason being that the term \(\cost_1(x_{k-1})-\costinf*\) therein is not necessarily positive (by the same argument, \(\mathcal L_k(x_\star)\) is not guaranteed to be positive).
			For this reason we talk in terms of \emph{quasi}-descent when referring to inequality \eqref{eq:puya:muk:descent}, a complication that, similarly to the analysis in \cite{solodov2007explicit,cabot2005proximal}, is the culprit of a nonstraightforward derivation of convergence results.
			Nevertheless, regardless of the sign of \(\cost_1(x_{k-1})-\costinf*\), the combination of \cref{ass:cost1}, boundedness of \(\seq{\alphk}\) and \(\seq{\rhok}\) (to be established later), the fact that \(x_k\in\dom\cost_2\) holds for all \(k\), and the slow control condition ensures that \(\seq{\mathcal L_k(x_\star)}\) is lower bounded (in fact, \(\liminf_{k\to\infty}\mathcal L_k(x_\star)\geq0\)) for any \(x_\star\in\X_1\).

		\subsection{Convergence recipe for proximal gradient iterations}\label{sec:recipe}

			The convergence of the two proposed algorithms hinges on the behavior of proximal gradient iterations when some implicit conditions are met.
			This is materialized through a convergence recipe relying on the following properties of the generated stepsize sequence.

			\begin{center}
				\setlength{\fboxsep}{8pt}%
				\setlength{\fboxrule}{2pt}%
				\fbox{\parbox{0.9\linewidth}{%
					\textsf{\bfseries Properties of stepsizes \(\alphk\) and inverse penalties \(\sigk\)}

					\medskip
					There exist \(\alpha_{\rm max}\geq\alpha_{\rm min}>0\) and \(\nu\in(0,1)\) such that, for every \(k\in\N\),

					\begin{minipage}[t]{0.57\linewidth}
						\vspace*{-.5\baselineskip}%
						\begin{alphaproperties}
						\item \label{cond:alphk}%
							\(\alphk*\leq\min\set{\balphk*,\alpha_{\rm max}}\) with \(\balphk*\) as in \eqref{eq:adabim:alphk*_0}%

						\item \label{cond:alphk:nu}%
							\(\alphk*\lk*\leq\nu\)

						\item \label{cond:alphamin}%
							\(\alphk*\geq\alpha_{\rm min}\)
						\end{alphaproperties}
					\end{minipage}
					\hfill\vline\hfill
					\begin{minipage}[t]{0.35\linewidth}
						\vspace*{-.5\baselineskip}%
						\begin{sigmaproperties}
						\item \label{cond:monotone}%
							\(0<\sigk* \leq \sigk\)

						\item \label{cond:slowcontrol}%
							\(\sigk\to 0\) and \(\sum_{k\in\N}\sigk=\infty\)
						\end{sigmaproperties}
					\end{minipage}
				}}%
			\end{center}

			Before presenting the unifying convergence recipe, we establish intermediate but crucial results such as boundedness of the sequence without imposing a uniform lower bound on the stepsize as in \cref{cond:alphamin}.
			Optimality of the limit points, however, will ultimately hinge on this final assumption and will be presented in \cref{thm:recipe}.

			\begin{lemma}\label{thm:PG}%
				Suppose that \cref{ass:basic} holds, and consider proximal gradient iterations \eqref{eq:PG} with \(\seq{\alphk}\) and \(\seq{\sigk}\) complying with \cref{cond:alphk,cond:alphk:nu,cond:monotone}.
				Then, the following hold:
				\begin{enumerate}
				\item \label{thm:PG:rhomax}%
					\(
						\rhok
					\coloneqq
						\frac{\sigk\alphk}{\sigma_{k-1}\alpha_{k-1}}
					\leq
						\rho_{\rm max}
					\coloneqq
						\max\set{
							\frac{\alpha_0}{\alpha_{-1}},\,
							\frac{1+\sqrt5}{2}
						}
					\)
					for every \(k\in\N\).

				\item \label{thm:PG:qSD}%
					For \(\mathcal L_k\) as in \cref{thm:LW} it holds that
					\(
						\mathcal L_{k+1}(x_\star)
					\leq
						\mathcal L_k(x_\star)
						-
						\sigk\alphk(1+ \rhok -\rhok*^2)
						\bigl(
							\cost_1(x_{k-1})-\costinf*
						\bigr)
					\)
					for all \(k\in\N\) and \(x_\star\in\X_1\).

				\item \label{thm:PG:barPhi}%
					\(
						\bphik*(x_{k+1})
					\leq
						\bphik(x_k)
						-
						\tfrac{1-\nu}{\alphk*}\|x_{k+1}-x_k\|^2
					\)
					holds for every \(k\in\N\).
					In particular, \(\|x_k-x_{k-1}\|\to0\) as \(k\to \infty\), \(\seq{\bphik(x_k)}\) is convergent, the following  worst-case rate holds
					\[
						\min_{k\leq K} \|x_{k+1} - x_k\|^2 \leq
						\tfrac{\alpha_{\rm max} \bar{\originalvarphi}_0(x_{0})}{(1-\nu)(1+K)},
					\]
					and both \(\seq{\sigk\cost_1(x_k)}\) and \(\seq{\cost_2(x_k)}\) are bounded.

				\item \label{thm:PG:boundedness}%
					The sequence \(\seq{x^k}\) is bounded.
				\end{enumerate}
			\end{lemma}
			\begin{proof}~
				\begin{proofitemize}
				\item \ref{thm:PG:rhomax}~
					Observing that
					\(
						\rhok*
					\leq
						\sqrt{\tfrac{\sigk}{\sigma_{k-1}}\bigl(1+\rhok\bigr)}
					\leq
						\sqrt{1+\rhok}
					\)
					by \cref{cond:alphk}, the assertion follows from a trivial induction argument.

				\item \ref{thm:PG:qSD}~
					Since \(\alphk\leq\balphk\), from the definition of \(\balphk\) in \eqref{eq:adabim:alphk*_0} it follows that the multiplying coefficients in \eqref{eq:puya:muk:descent} are positive:
					\begin{align*}
						\mathcal L_{k+1}(x_\star) - \mathcal L_k(x_\star)
					\leq{} &
						-
						\sigk*\alphk*
						\left(
							\overbracket*[0.5pt]{
								\tfrac{\Delk}{\rhok*}
							}^{\geq0}
							+
							\tfrac{1}{\sigk}
							\bigl(
								\overbracket*[0.5pt]{
									1+\rhok*
									-
									\rho_{k+2}^2\tfrac{\sigk}{\sigk*}
								}^{\geq0}
							\bigr)
						\right)
						\bcost_2(x_k)
					\\
					&
						-
						\sigk\alphk
						\bigl(
							\overbracket[0.5pt]{
								1+\rhok-\rhok*^2
							}^{\geq0}
						\bigr)
						\bigl(\cost_1(x_{k-1})-\costinf*\bigr)
					\\
					&
						-
						\left(
							\underbracket*[0.5pt]{
								\tfrac{1}{4}
								+
								\rhok*^2(\alphk\lk-\alphk^2\Lk^2)
								-
								\alphk\sigk\lk^{(2)}\Delk
							}_{\geq0}
						\right)
						\|x_{k-1}-x_k\|^2.
						\vphantom{
							\underbracket[0.5pt]{
								\tfrac{1}{4}
							}_0
						}
					\end{align*}
					The claim then readily follows from the fact that \(\bcost_2\geq0\), cf. \cref{thm:bar>=0}.

				\item \ref{thm:PG:barPhi}~
					We have
					\begin{align*}
						\bphik*(x_{k+1})
						-
						\bphik(x_k)
					\leq{} &
						\bphik*(x_{k+1})
						-
						\bphik*(x_k)
					\\
					={} &
						\phik*(x_{k+1})
						-
						\phik*(x_k)
					\\
					={} &
						\sigk*\left(
							\Phik*(x_{k+1})
							-
							\Phik*(x_k)
						\right)
					\overrel*[\leq]{\eqref{eq:phikdescent}}
						- \tfrac{1-\nu}{\alphk*}
						\|x_{k+1}-x_k\|^2.
						\numberthis \label{eq:descent:bphik}
					\end{align*}
					Here, the first inequality uses the fact that \(\sigk\geq\sigk*\), and therefore \(\bphik=\sigk\bcost_1+\bcost_2\) is smaller than \(\bphik*=\sigk*\bcost_1+\bcost_2\) on \(\dom\cost_2\), cf. \cref{thm:bar>=0}, and \(x_k\in\dom\gk=\dom\cost_1\cap\dom\cost_2\subseteq\dom\cost_2\);
					the second inequality follows from the fact that \(\alphk*\lk*\leq\nu\).
					This shows the sought inequality.
					In turn, since \(\nu<1\), the sequence \(\seq{\bphik(x_k)}\) is decreasing; since \(x_k\in\dom\cost_2\), it follows from \cref{thm:bar>=0} that \(\bphik(x_k)=\sigk\bcost_1(x^k)+\bcost_2(x^k)\geq0\), hence that it is convergent and that the positive-valued sequences \(\seq{\sigk\bcost_1(x^k)}\) and \(\seq{\bcost_2(x^k)}\) are bounded.

					The convergence rate result follows immediately by telescoping \eqref{eq:descent:bphik}
					\begin{equation}\label{eq:sublin:res}
						\tfrac{(1-\nu)(K+1)}{\alpha_{\rm max}}\min_{k\leq K}\|x_{k+1}-x_k\|^2
					\leq
						\sum_{k=0}^{K}\tfrac{1-\nu}{\alphk*}\|x_{k+1}-x_k\|^2
					\leq
						\bar{\originalvarphi}_0(x_{0}),
					\end{equation}
					where \cref{cond:alphk} was used to bound \(\alphk*\leq\alpha_{\rm max}\).

				\item \ref{thm:PG:boundedness}~
					We pattern the proof structure of \cite[Thm. 3.2]{solodov2007explicit}, thereby considering two mutually exclusive cases.
					\begin{proofitemize}[leftmargin=*,topsep=0pt]
					\item
						{\it Case 1: \(\cost_1(x_k)\geq\costinf*\) holds for \(k\) large enough.}

						Recall the definition of \(W_k\) and \(\mathcal L_k\) in \cref{thm:LW}.
						Observing that \(W_k\geq0\) and \(\frac12\|x_k-x_\star\|^2\leq\mathcal L_k(x_\star)\), \cref{thm:PG:qSD} implies that \(\seq{\mathcal L_k(x_\star)}\) converges and that consequently \(\seq{x_k}\) is bounded.

					\item
						{\it Case 2: \(\cost_1(x_k)<\costinf*\) holds infinitely often.}

						In this case, for every \(k\) large enough the index
						\begin{equation}\label{eq:ik}
							i_k
						\coloneqq
							\max\set{i\leq k}[\cost_1(x_i)<\costinf*]
						\end{equation}
						is well defined.
						We proceed by intermediate claims.
						\renewcommand{\claimscounter}{\ref*{thm:PG}}%
						\begin{claims}
						\item\label{claim:convergence:ik}%
							{\it the sequences \(\seq{x_{i_k}}\) and \(\seq{x_{i_k+1}}\) are bounded}.

							The optimal set
							\[
								\X_1
							=
								\set{x\in\X_2}[\cost_1(x)\leq \costinf*]
							=
								\set{x}[\cost_1(x)\leq\costinf*,\ \bcost_2(x)\leq0]
							\]
							coincides with a sublevel set of the convex function \(h\coloneqq\max\set{\cost_1-\costinf*,\bcost_2}\).
							Since it is nonempty and bounded by \cref{ass:X1}, \(h\) is level bounded; see, e.g., \cite[Prop. 11.13]{bauschke2017convex} or \cite[Lem. 1]{themelis2019acceleration}.
							Note that \cref{thm:PG:barPhi} implies that \(\seq{\bcost_2(x_k)}\) is (upper) bounded, which combined with the fact that \(\cost_1(x_{i_k})<\costinf*\) implies that \(\seq{x_{i_k}}\) lies in a sublevel set of \(h\), and is therefore bounded.
							In turn, \cref{thm:PG:barPhi} implies that so is \(\seq{x_{i_k+1}}\).
						\item
							{\it the whole sequence \(\seq{x_k}\) is bounded}.

							To this end, it remains to show that \(\seq{x_{k'}}[k'\in K']\) is bounded, where
							\begin{equation}\label{eq:K'}
								K'\coloneqq\set{k'\in\N}[k'\geq i_{k'}+2].
							\end{equation}
							With \(\rho_{\rm max}\) as in \cref{thm:PG:rhomax}, for every \(k\in\N\)
							\begin{align*}
								\mathcal L_{i_k+1}(x_\star)
							\leq{} &
								\overbracket[0.5pt]{
									\tfrac12
									\|x_{i_k+1}-x_\star\|^2
									\vphantom{\left(\tfrac{\sigma_{i_k+1}}{\sigma_{i_k}}\right)}
								}^{\text{bounded by \cref{claim:convergence:ik}}}
								+~~
								\overbracket*[0.5pt]{
									\tfrac{1}{4}\|x_{i_k+1}-x_{i_k}\|^2
									\vphantom{\left(\tfrac{\sigma_{i_k+1}}{\sigma_{i_k}}\right)}
								}^{\mathclap{\to0 \text{ by \cref{thm:PG:barPhi}}}}
								+
								\overbracket[0.5pt]{
									\sigma_{i_k+1}\alpha_{i_k+1}(1+\rho_{i_k+1})
									(\cost_1(x_{i_k})-\costinf*)
									\vphantom{\left(\tfrac{\sigma_{i_k+1}}{\sigma_{i_k}}\right)}
								}^{<0}
							\\
							&
								+
								\underbracket[0.5pt]{
									\alpha_{\rm max}\rho_{\rm max}^2
									\bcost_2(x_{i_k})
									+
									\alpha_{\rm max}\left(1-\tfrac{\sigma_{i_k+1}}{\sigma_{i_k}}\right)
									\bcost_2(x_{i_k+1})
								}_{\text{bounded by \cref{thm:PG:barPhi}}}.
							\end{align*}
							In particular, we have
							\begin{equation}\label{eq:supL}
								\sup_{k\in\N}\mathcal L_{i_k+1}(x_\star)<\infty.
							\end{equation}
							Let now \(k'\in K'\).
							Since \(\cost_1(x_j)-\costinf*\geq0\) for \(j=i_{k'}+1,\dots,k'\), observe that
							\begin{equation}\label{eq:Lk'geq}
								\tfrac12\|x_{k'}-x_\star\|^2
							\leq
								\mathcal L_{k'}(x_\star)
							\quad
								\forall k'\in K',
							\end{equation}
							and \cref{thm:PG:qSD} yields that
							\begin{align*}
								\mathcal L_{k'}(x_\star)
							\leq{} &
								\mathcal L_{k'-1}(x_\star)
							\leq\dots\leq
								\mathcal L_{i_{k'}+2}(x_\star)
							\numberthis\label{eq:Lk'descent}
							\\
							\leq{} &
								\mathcal L_{i_{k'}+1}(x_\star)
								-
								\sigma_{i_{k'}+1}\alpha_{i_{k'}+1}(1+\rho_{i_{k'}+1}-\rho_{i_{k'}+2}^2)
								\bigl(
									\cost_1(x_{i_{k'}})-\costinf*
								\bigr)
							\\
							\leq{} &
								\mathcal L_{i_{k'}+1}(x_\star)
								+
								\underbracket[0.5pt]{
									\tfrac{\sigma_{i_{k'}+1}}{\sigma_{i_{k'}}}
								}_{\leq1}
								\underbracket[0.5pt]{
									\vphantom{\tfrac{\sigma_{i_{k'}+1}}{\sigma_{i_{k'}}}}
									\alpha_{\rm max}(1+\rho_{i_{k'}+1}-\rho_{i_{k'}+2}^2)
								}_{\text{\clap{bounded by \cref{thm:PG:rhomax}}}}
								\underbracket[0.5pt]{
									\vphantom{\tfrac{\sigma_{i_{k'}+1}}{\sigma_{i_{k'}}}}
									\sigma_{i_{k'}}
									\bigl|
										\cost_1(x_{i_{k'}})-\costinf*
									\bigr|
								}_{\text{\clap{bounded by \cref{thm:PG:barPhi}}}}
								< \infty
							\qquad
								\forall k'\in K',
							\end{align*}
							where \eqref{eq:supL} was used in the last inequality.
							Here, boundedness of the under-bracketed term follows from boundedness of \(\seq{x_{i_{k'}}}[k'\in K']\) and lower semicontinuity of \(\cost_1\).
							Then, \eqref{eq:Lk'geq} implies that the sequence \(\seq{x_{k'}}[k'\in K']\) is bounded.
							Combined with \cref{claim:convergence:ik} and the fact that the index set \(K'\) is the complement of the indices therein, the claim follows.
						\qedhere
						\end{claims}
					\end{proofitemize}
				\end{proofitemize}
			\end{proof}

			\begin{theorem}[convergence recipe for proximal gradient iterations]\label{thm:recipe}%
				Suppose that \cref{ass:basic} holds, and consider the proximal gradient iterations \eqref{eq:PG} with \(\seq{\sigk}\) and \(\seq{\alphk}\) complying with all \cref{cond:alphk,cond:alphk:nu,cond:alphamin,cond:monotone,cond:slowcontrol}.
				Then,
				\begin{enumerate}
				\item \label{thm:recipe:solution}%
					 (\(\seq{x_k}\) is bounded and) \(\seq{\dist(x_k, \X_1)}\) converges to zero.
				\item \label{thm:recipe:rate}%
					both \(\bphik(x_k)\) and \(W_k\) as in \cref{thm:LW} converge to 0 as \(k\to\infty\), and the following worst-case rate holds
					\[
						\min_{k\leq K}\dist^2(0, \partial \phik*(x_{k+1})
					\leq
						\tfrac{ (1+ \alpha_{\rm max}L_{f_0, \V})^2}{\alpha_{\rm min}(1-\nu)(K+1)}\bar{\originalvarphi}_0(x_0),
					\]
					where \(L_{f_0, \V}\) is a Lipschitz modulus for \(\sigma_0\nabla\f_1+\nabla\f_2\) on \(\V\coloneqq\conv\set{x^k}[k\in\N]\).
				\end{enumerate}
			\end{theorem}
			\begin{proof}
				We begin by remarking that boundedness of the sequence \(\seq{x_k}\) is ensured by \cref{cond:alphk,cond:alphk:nu}, as shown in \cref{thm:PG:boundedness}.
				We next prove each claim individually.
				\begin{proofitemize}
				\item \ref{thm:recipe:rate}~
					Consider a convergent subsequence \(x_{k_j}\to x_\infty\), so that \(x_{k_j+1}\to x_\infty\) by \cref{thm:PG:barPhi}.
					Up to further extracting if necessary we have that
					\(\alpha_{k_j+1}\to\alpha_\infty\geq\alpha_{\rm min}>0\) and \(\sigma_{k_j}\to0\).
					Observe that
					\begin{align*}
						x_{k+1}
					={} &
						\argmin h({}\cdot{};x_k,\alpha_{k+1},\sigk*),
					\shortintertext{where}
						h(w;x,\alpha,\sigma)
					\coloneqq{} &
						(\sigma \g_1+\g_2)(w)
						+
						\tfrac{1}{2\alpha}
						\bigl\|w-x+\alpha\nabla\bigl(\sigma \f_1+\f_2\bigr)(x)\bigr\|^2
					\end{align*}
					is level bounded in \(w\) locally uniformly in \((x,\alpha,\sigma)\), as a function from \(\R^n\times\bigl(\R^n\times[\alpha_{\rm min},\infty)\times[0,\infty)\bigr)\) to \(\Rinf\).
					Since \(h\) is continuous in \((x,\alpha,\sigma)\), it follows from \cite[Thm. 1.17]{rockafellar2011variational} that
					\[
						x_\infty
					\in
						\argmin h({}\cdot{};x_\infty,\alpha_\infty,0)
					=
						\prox_{\alpha_\infty \g_2}\bigl(x_\infty-\alpha_\infty\nabla \f_2(x_\infty)\bigr),
					\]
					this condition being equivalent to \(x_\infty\in\argmin(\f_2+\g_2)\defeq\X_2\).
					We next show that \(\bphik(x_k)\to0\).

					For any \(x_\star\in\X_1\), the subdifferential characterization of \(x_{k+1}\in\prox_{\alphk*\gk*}(x_k-\alphk*\nabla\fk*(x_k))\)
					\[
						\tfrac{x_{k-1}-x_k}{\alphk}
						-
						(\nabla\fk(x_{k-1})-\nabla\fk(x_k))
					\in
						\partial\phik(x_k)
					\]
					implies that
					\[
						\phik(x_\star)
					\geq
						\phik(x_k)
						+
						\innprod{
							\overbracket*[0.5pt]{
								\tfrac{x_{k-1}-x_k}{\alphk}
								-
								(\nabla\fk(x_{k-1})-\nabla\fk(x_k))
							}^{\to0 \text{ by \cref{thm:PG:barPhi}}}
						}{
							\overbracket*[0.5pt]{
								x_\star-x_k
								\vphantom{\tfrac{x_{k-1}-x_k}{\alphk}}
							}^{\text{\clap{bounded}}}
						},
					\]
					hence that
					\[
						\costinf_2
					=
						\lim_{k\to\infty}\phik(x_\star)
					\geq
						\limsup_{k\to\infty}\phik(x_k)
					=
						\limsup_{k\to\infty}\left(\sigk\bcost_1(x_k)+\overbracket*[0.5pt]{\cost_2(x_k)}^{\geq\costinf_2}\right)
					\geq
						\limsup_{k\to\infty}\sigk\bcost_1(x_k)
						+
						\costinf_2.
					\]
					Here, the second equality is obtained by adding and subtracting \(\sigk\costinf_1\) along with the fact that \(\sigk\to0\).
					Since \(\sigk\bcost_1\geq0\), necessarily \(\sigk\bcost_1(x_k)\to0\) and consequently \(\limsup_{k\to\infty}\cost_2(x_k)=\costinf_2\).
					Similarly, since \(\cost_2\geq\costinf_2\) this necessarily implies that \(\cost_2(x_k)\to\costinf_2\).
					This shows that \(\bcost_2(x^k)\to0\) as claimed.

					We next show that \(\seq{W_k}\) too is vanishing.
					Observe that
					\begin{align*}
						W_k
					\defeq{} &
						\tfrac{1-4\alphk\sigk\lk^{(2)}\Delk}{4}
						\overbracket[0.5pt]{
							\|x_k-x_{k-1}\|^2
						}^{\mathclap{\to 0\text{ by \cref{thm:PG:barPhi}}}}
						+
						\tfrac{\sigk*}{\sigk}\alphk*\rhok*
						\overbracket[0.5pt]{
							\bcost_2(x_{k-1})
						}^{\mathclap{\to 0}}
					\\
					&
						+
						\underbracket[0.5pt]{
							\sigk\alphk\Delk\bcost_2(x_k)
						}_{\to 0}
						+
						\alphk
						\left(1+\rhok\right)
						\underbracket[0.5pt]{
							\sigk
							\bigl(
								\cost_1(x_{k-1})
								-
								\costinf*
							\bigr)
						}_{\to0}
					\to
						0
					\quad\text{as \(k\to\infty\)},
					\end{align*}
					where we also used the fact that \(\seq{\alphk}\leq\alpha_{\rm max}\) and \(\seq{\rhok}\leq\rho_{\rm max}\) as in \cref{thm:PG:rhomax}.

					The asserted sublinear rate is a consequence of \eqref{eq:sublin:res}, the characterization of the residual in \cref{thm:H} (see \eqref{eq:Hk:phi}), and the upper bounds for \(\lk\), \(\Lk\) in \cref{thm:l<=L}:
					\begin{align*}
						\textstyle\sum_{k=0}^K\dist^2(0, \partial \phik*(x_{k+1})
					\leq{}&
						\textstyle\sum_{k=0}^K\tfrac{1}{\alphk*^2}\|\Hk*(x_{k})-\Hk*(x_{k+1}) \|^2
					\\
					\leq {}&
						\tfrac{M_{\rm max}^2}{\alpha_{\rm min}}\textstyle\sum_{k=0}^K \tfrac{1}{\alphk*}\|x_{k}-x_{k+1}\|^2
					\leq
						\tfrac{M_{\rm max}^2}{\alpha_{\rm min}(1-\nu)}\bar{\originalvarphi}_0(x_{0}),
					\end{align*}
					where \(M_{\rm max} = \sup_{k\in\N}M_k\leq 1+ \alpha_{\rm max}L_{f_0, \V}\), see \cref{thm:lL}.

				\item \ref{thm:recipe:solution}~
					We begin by observing that, since \(\bcost_2(x^k)\to0\), by lower semicontinuity all limit points of \(\seq{x^k}\) belong to \(\X_2\).
					As done in the proof of \cref{thm:PG:boundedness} we consider two possible cases.
					\begin{proofitemize}[leftmargin=*,topsep=0pt]
					\item
						{\it Case 1: \(\cost_1(x_k)\geq\costinf*\) holds for \(k\) large enough.}

						In this case, we will actually show that the sequence \(\seq{x_k}\) converges to a solution of \eqref{eq:P}.
						We first argue that there exists an optimal limit point; to this end, since all limit points are feasible and because of lower semicontinuity, it suffices to show that \(\liminf_{k\to\infty}\cost_1(x_k)=\costinf*\).
						A telescoping argument on \cref{thm:PG:qSD} along with \cref{cond:alphamin} yields
						\begin{equation}
							\alpha_{\rm min}\sum_{k\in \N}\sigk(1+\rhok-\rhok*^2)\bigl(\cost_1(x_{k-1})-\costinf*\bigr)
						\leq
							\sum_{k\in \N}\sigk\alphk(1+\rhok-\rhok*^2)\bigl(\cost_1(x_{k-1})-\costinf*\bigr)
						<
							\infty.
						\end{equation}
						Since \(\sum_{k\in\N}\sigk=\infty\), necessarily
						\begin{align*}
							\liminf_{k\to\infty}(1+\rhok-\rhok*^2)(\cost_1(x_{k-1})-\costinf*)
						=
							0.
							\numberthis\label{eq:liminf:cost1}
						\end{align*}
						Notice that \(\limsup_{k\to\infty}(1+\rhok-\rhok*^2)>0\), for otherwise \(1+\rhok-\rhok*^2 \to 0\), implying that \(\liminf_{k \to \infty} \rhok >0\), and consequently that \(\alphk\sigk\) eventually increases exponentially, which contradicts \(\sigk\alphk \leq \sigk\alpha_{\rm max} \to 0\).
						Hence, \(\limsup_{k\to\infty}(1+\rhok-\rhok*^2)>0\), which along with \eqref{eq:liminf:cost1} implies that \(\liminf_{k\to\infty}\cost_1(x_{k-1})=\costinf*\).
						Therefore, an optimal limit point exists, be it \(x_\infty\).
						Since \(W_k\to0\) by assertion \ref{thm:recipe:rate} and since \(\mathcal L_k(x_\infty)\) converges, it follows that \(\frac12\|x_k-x_\infty\|^2 = \mathcal L_k(x_\infty) - W_k\) converges as well.
						Since along a subsequence \(\frac12\|x_k-x_\infty\|^2\) converges to zero, necessarily \(\lim_{k\to\infty}\frac12\|x_k-x_\infty\|^2=0\), proving that the entire sequence \(\seq{x_k}\) converges to \(x_\infty\).
					\item
						{\it Case 2: \(\cost_1(x_k)<\costinf*\) holds infinitely often.}

						Recall the index \(i_k\) and the set \(K'\) defined in \eqref{eq:ik} and \eqref{eq:K'}.
						Having established boundedness of the entire sequence, limit points of \(\seq{x_{i_k}}\) exist and, as shown above, belong to \(\X_2\).
						Moreover, it follows from \eqref{eq:ik} that \(\limsup_{k\to\infty}\cost_1(x_{i_k})\leq\costinf*\), and a lower semicontinuity argument then yields that all the limit points of \(\seq{x_{i_k}}\) attain the optimal cost and are therefore optimal.
						Since \(x_{i_k}-x_{i_k+1}\to0\) by \cref{thm:PG:barPhi}, the same is also true for \(\seq{x_{i_k+1}}\), and in particular
						\(\dist(x_{i_{k}+1},\X_1) \to 0\) as \(k \to \infty\).

						For each \(k\) let \(\bar x_k\coloneqq\proj_{\X_1}\!x_k\), which is well defined since \(\X_1\neq\emptyset\) is closed and convex.
						Recalling that \(\mathcal L_k(x_\star)=W_k + \tfrac12\|x_k-x_\star\|^2\), for \(k'\in K'\) we have
						\begin{align*}
							\tfrac12\dist(x_{k'},\X_1)^2
						\leq{} &
							\tfrac12\|x_{k'}-\bar x_{i_{k'}+1}\|^2
						={}
							\mathcal L_{k'}(\bar x_{i_{k'}+1})-W_{k'}
						\\
							\dueto{\scriptsize by \eqref{eq:Lk'descent}}
							\
						\leq{} &
							\mathcal L_{i_{k'}+2}(\bar x_{i_{k'}+1})
							-
							W_{k'}
						={}
							\tfrac12\|x_{i_{k'}+2} - \bar x_{i_{k'}+1}\|^2
							+
							W_{i_{k'}+2}
							-
							W_{k'}
						\\
						\leq{} &
							\tfrac12\big(\dist(x_{i_{k'}+1},\X_1) + \|x_{i_{k'}+1} - x_{i_{k'}+2}\|\big)^2
							+
							W_{i_{k'}+2}
							-
							W_{k'}
						\to
							0
						\end{align*}
						as \(K'\ni k'\to\infty\), where the limit follows from \cref{thm:PG:barPhi}, the vanishing of \(\seq{W_k}\) established in assertion \ref{thm:recipe:rate}, and the fact that \(i_k\to\infty\) as \(k\to\infty\).
						Since \(K'\cup\set{i_k,i_k+1}[k\in\N]=\N\), we conclude that all the limit points of \(\seq{x^k}\) are optimal.
					\qedhere
					\end{proofitemize}
				\end{proofitemize}
			\end{proof}

	\section{Simulations}\label{sec:num}

		In this section the performance of the proposed algorithms is evaluated through a series of simulations on standard problems on both synthetic data and standard datasets from the LIBSVM dataset \cite{chang2011libsvm}.
		All the algorithms are implemented in the Julia programming language and are available online.\footnote{%
			\url{https://github.com/pylat/adaptive-bilevel-optimization}
		}
		An overview of the algorithms included in the simulations is provided in the following subsection.

		In accounting for the difference in iteration complexity among the methods, the simulations report the progress against the number of calls to \(\nabla\f_2\), since in all problems \(\nabla\f_1\) and proximal operations have negligible cost.
		As explained in \cref{sec:solodov}, this criterion favors the method \solodov{}, as it ignores the cost of the backtracks which involve function evaluations.
		On the contrary, all the backtracking steps included in \refadabim*, which involve gradient evaluations, are fully accounted for in the comparisons.

		\subsection{Compared algorithms}

			When applicable, other than \refadabim* and \refstabim* the algorithms involved in the simulations are \solodov, \bigsam, and \itthreeD.
			For \itthreeD{}, \(\sigk=\nicefrac{1}{(k+1)^2}\) was used,{}  for \bisg{}, \(\sigk = \nicefrac{1}{(k+1)^p}\) with \(p=0.95\), and \(\sigk=\nicefrac{1}{k+1}\) was adopted for the rest.
			Although only \bigsam{} and \bisg{} require \(\sigk\in(0,1]\), this limitation was applied across all methods to maintain more uniform comparisons.

			It is also worth noting that \refadabim* is not sensitive to the choice of initial stepsize \(\alpha_0\), as future values are automatically adjusted during the iterations.
			In all the simulations, the parameter \(\nu = 0.99\) was used, and
			the parameter \(\alpha_{\rm max}\) appearing in \eqref{eq:adabim:alphk*_0} of \refadabim* was set as a large constant; as remarked before, it is only of theoretical significance.
			As we will see, the same parameter in \solodov{} is instead crucial for dictating the algorithmic performance.
			These facts are better detailed in the following brief description of the algorithms compared against in the simulations.

			\subsubsection{Solodov's explicit descent method (\solodov-\texorpdfstring{\(r\)}{r})}\label{sec:solodov}%
				This is Solodov's explicit descent method \cite[Alg. 2.1]{solodov2007explicit} already outlined in \eqref{eq:solodov}.
				In the simulations, the suffix ``-\(r\)'' is used to distinguish different choices for the value of \(\widehat\alpha_0\) therein; namely, whenever \(\f^2\) is \(L_{\f_2}\)-Lipschitz differentiable, we set \(\widehat\alpha_0\coloneqq\frac{r}{L_{\f_2}}\).
				In addition to \cref{ass:basic}, the algorithm requires:
				\begin{itemize}[itemsep=0pt]
				\item
					\(\g_1=0\);
				\item
					\(\g_2=\indicator_D\) for a nonempty, closed, and convex set \(D\subseteq\R^n\).
				\end{itemize}
				The backtracks involved in \solodov{} do not require additional gradient evaluations, but instead require function evaluations which are not reflected in the comparisons in terms of total number of gradients.
				For this reason, in the top row of \cref{fig:gam} we provided a sample plot for three selected applications demonstrating the higher number of backtracks that it incurs compared to \refadabim*.
				As evident in the figures, in practice \solodov{} is sensitive to parameter tuning; while selecting a larger \(\widehat\alpha_0\) can lead to larger stepsizes and consequently faster convergence speed in terms of number of iterations (gradient evaluations), it results in a higher number of backtracks (each requiring one cost evaluation), and vice versa.

				In all the simulations, both for \solodov{} and \refadabim* we used the backtrack parameter \(\eta = \nicefrac12\), and the linesearch related parameter \(\nu=0.99\).

			\subsubsection{Bilevel gradient sequential averaging method (\bigsam)}\label{sec:bigsam}
				Proposed in \cite{sabach2017first}, \bigsam{} addresses strongly convex bilevel methods under global Lipschitz differentiability assumptions.
				Specifically, in addition to \cref{ass:basic} the algorithm requires that
				\begin{itemize}[itemsep=0pt]
				\item
					\(\f_1\) is \(L_{\f_1}\)-Lipschitz differentiable and \(\mu_{\f_1}\)-strongly convex;
				\item
					\(\f_2\) is \(L_{\f_2}\)-Lipschitz differentiable;
				\item
					\(\g_1=0\).
				\end{itemize}
				\Bigsam{} iterates
				\[
					\begin{cases}
						\xk*_1 & {}= \xk_1-\alpha^{(1)}\nabla\f_1\left(\xk_1\right)\\
						\xk*_2 & {}= \prox_{\alpha^{(2)}\g_2}\left(\xk_1-\alpha^{(2)}\nabla\f_2\left(\xk_2\right)\right)\\
						x_{k+1}& {}= \sigk*\xk*_1+(1-\sigk*)\xk*_2,
					\end{cases}
				\]
				where \(\alpha^{(1)}\leq\frac{2}{L_{\f_1}+\mu_{\f_1}}\), \(\alpha^{(2)}\leq\frac{1}{L_{\f_2}}\).
				Although in \bigsam{} the sequence \(\seq{\sigk}\) has a different interpretation that the one in \eqref{eq:phik}, it must still comply with \eqref{eq:sigk} and in addition \(\sigma_1\leq1\).
				For this reason we opted to use the same notation.

			\subsubsection{\ItthreeD}\label{sec:itthreeD}
				This method, presented in \cite{garrigos2018iterative}, is specialized to linear inverse problems and operating on the dual formulation, which complicates the comparison with the other methods.
				A strongly convex upper layer cost is required, but the iterations only involve gradient (and not proximal) evaluations on its (Lipschitz-differentiable) conjugate.
				Similarly, the lower level cost is an infimal convolution between a prox-friendly and a strongly convex function, making its dual the sum of a prox-friendly and a Lipschitz-differentiable terms.
				The requirements on the primal formulation are roughly as follows:
				\begin{itemize}[itemsep=0pt]
				\item
					\(\cost_1\) is \(\mu_{\cost_1}\)-strongly convex;
				\item
					\(\cost_2\) is a coercive ``data-fit'' function.
				\end{itemize}
				In referring the reader to \cite{garrigos2018iterative} for a rigorous account on the problem formulation and its requirements, we point out that among our simulations \itthreeD{} is only applicable to the linear inverse problem of \cref{sec:linverse} with \(\ell^2\)-norm upper layer cost.
				In that setting, initializing with \(x_0\in\range\trans A\) the method performs the following iterations
				\[
					x_{k+1}
				=
					x_k
					-
					\gamma\nabla\fk(x_k)
				=
					x_k
					-
					\trans A(Ax_k - b) - \gamma\sigk x_k,
				\]
				 each increasing the total gradient count by one.
				Remarkably, in this setting it does not constrain \(\seq{\sigk}\) to a nonsummable decay; see \cite[Rem. 10]{garrigos2018iterative}.
				For this reason, in the simulations inverse penalties \(\sigk=\nicefrac{1}{(k+1)^2}\) were used for \itthreeD, while \(\sigk=\nicefrac{1}{k+1}\) for all other methods.

			\subsubsection{Bi-Sub-Gradient method (\bisg)}\label{sec:bisg}
				Proposed in \cite{merchav2023convex}, this method can also cope with nondifferentiable terms in the upper level.
				It comes in two versions, depending on the type of operations on the upper level; we here consider the second one, as it relies on milder assumptions and is compatible with our proximal-gradient setting (the first one involves subgradient operations on the upper level).
				The standing assumptions are the following:
				\begin{itemize}[itemsep=0pt]
				\item
					\(\f_1\) is \(L_{\f_1}\)-Lipschitz differentiable;
				\item
					\(\f_2\) is \(L_{\f_2}\)-Lipschitz differentiable;
				\item
					\(\cost_1\) is coercive.
				\end{itemize}
				(An additional technical assumption of real valuedness of \(\cost_1\) is also imposed which nevertheless does not cause any loss of generality.)
				Notice that both our \cref{ass:cost1,ass:X1} are implied by (the existence of solutions and) the coercivity requirement on \(\cost_1\).
				The method alternates proximal-gradient operations with constant stepsize:
				\[
					\begin{cases}
						y_{k+1} & {}= \prox_{\alpha^{(2)}\g_2}\left(x_k-\alpha^{(2)}\nabla\f_2(x_k)\right)
					\\[3pt]
						x_{k+1} & {}= \prox_{\sigk*\alpha^{(1)}\g_1}\left(y_{k+1}-\sigk*\alpha^{(1)}\nabla\f_1(y_{k+1})\right),
					\end{cases}
				\]
				with \(\alpha^{(1)}=\frac{1}{\max\set*{1,L_{\f_1}}}\), \(\alpha^{(2)}=\frac{1}{L_{\f_2}}\), and \(\sigk*=\frac{1}{(k+1)^p}\) for some \(p\in(\frac12,1)\).
				As suggested in \cite{merchav2023convex}, \(p = 0.95\) was used in all the simulations.

		\subsection{Numerical experiments}\label{sec:numsim}

			We compare the algorithms on three benchmark bilevel problems with Lipschitz differentiable and strongly convex upper level cost \(\cost_1\), so as to satisfy the requirements of all the algorithms being compared.
			For two of these we also consider minimum \(\ell^1\)-norm versions in which the upper level \(\cost_1=\|{}\cdot{}\|_1\) is neither smooth nor strongly convex (but is nevertheless coercive).
			For these latter ones, only \refadabim*, \refstabim* and \bisg{} are applicable.

			\subsubsection{Logistic regression}\label{sec:logreg}
\begin{figure}[ht]
					\begin{subfigure}[b]{0.745\linewidth}
						\includetikz[width=\linewidth]{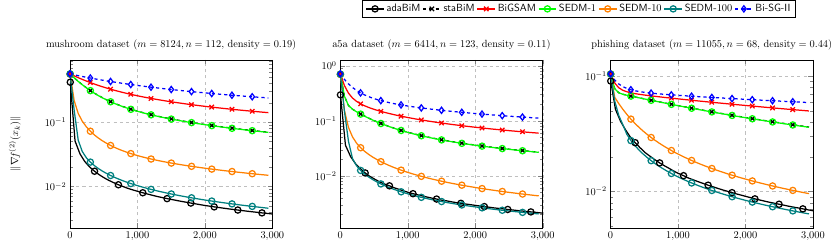}%
						\includetikz[width=\linewidth]{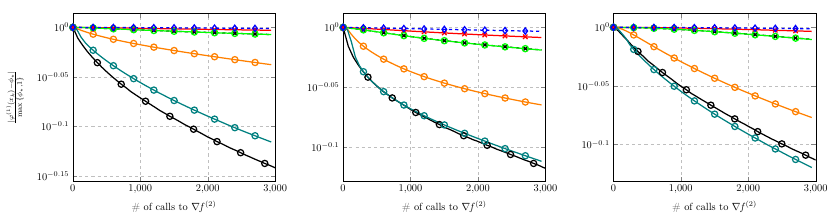}%
						\caption[]{\(\cost_1 = \tfrac12\|\cdot\|^2\)}
						\label{fig:logregL2}%
					\end{subfigure}
					\hfill
					\begin{subfigure}[b]{0.245\linewidth}
						\includetikz[width=\linewidth]{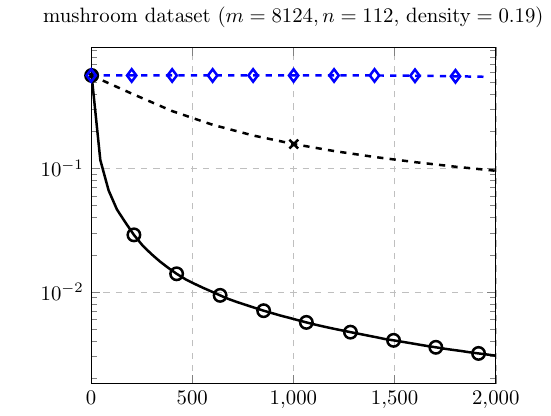}%
						\includetikz[width=\linewidth]{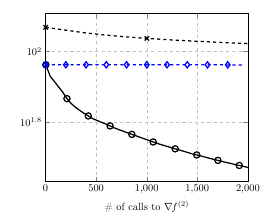}%
						\caption[]{\(\cost_1 = \|\cdot\|_1\)}
						\label{fig:logregL1}%
					\end{subfigure}
					\caption[]{Logistic regression problems of \cref{sec:logreg} with minimum \(\ell^p\)-norm solution, \(p=1,2\).}
				\end{figure}

				We consider the logistic regression problem
				\begin{subequations}\label{eq:logreg}
					\begin{align}
						\minimize_{x\in\R^n}{} &~ \cost_1(x)
					\\
						\stt{} &~ x \in \argmin_{w\in\R^n}
						\set{%
							\tfrac1m \textstyle\sum_{i=1}^m \left(y_i \log(s_i(w)) + (1-y_i) \log(1-s_i(w))\right)
							},
					\end{align}
				\end{subequations}
				where \(m,n\) are the number of samples and features, the pair \(a_i\in\R^{n+1}\) denotes the \(i\)-th sample (up to absorbing the bias terms), \(y_i\in\set{-1,1}\) is the associated label, and \(s_i(x) = (1+ \exp(-\trans{a_i}x))^{-1}\) is the logistic sigmoid function.
				In the simulations we used \(\cost_1 = \tfrac12\|{}\cdot{}\|^2\) (\cref{fig:logregL2}) and \(\cost_1 = \|{}\cdot{}\|_1\) (\cref{fig:logregL1}); for the latter, only one dataset is reported, as the plots for other ones are very similar.
				Note also that in the simulations for \refadabim*, \refstabim*{}  and \bisg{} we set \(\g_1 = \cost_1\) and \(\f_1 \equiv0\).
				For other methods (applicable only when \(\cost_1=\tfrac12\|{}\cdot{}\|^2\)) the upper level cost is captured using \(\f_1\) with \(L_{\f_1} = 1\) and its strong convexity modulus equal to \(1\) (in the case of \bigsam{}).
				For methods that require Lipschitz modulus of \(\nabla \f_2\), \(L_{\f_2} = \tfrac1{4m}\|A\|^2\) was used where \(A\) is the data matrix that is the concatenation of \((a_j)_{j=1,\dots,m}\).

			\subsubsection{Linear inverse problems with simulated data}\label{sec:linverse}
\begin{figure}[ht]
					\begin{subfigure}[b]{0.75\linewidth}
						\includetikz[width=\linewidth]{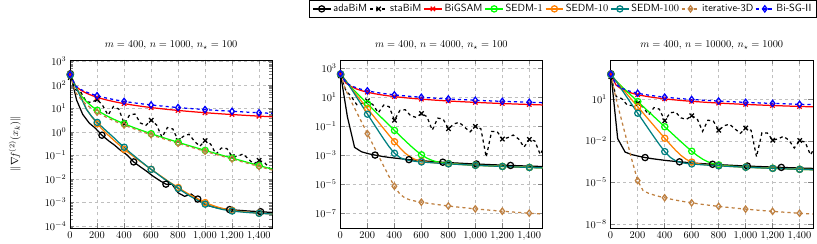}
						\includetikz[width=\linewidth]{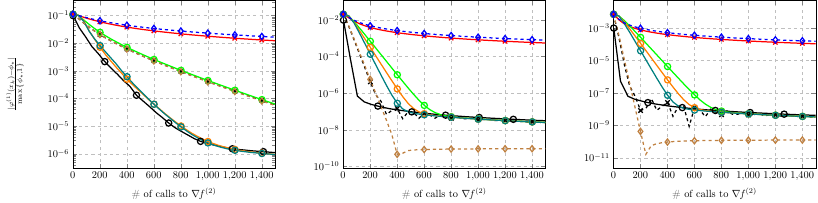}%
						\caption[]{\(\cost_1 = \tfrac12\|\cdot\|^2\)}
						\label{fig:linverseL2}%
					\end{subfigure}
					\hfill
					\begin{subfigure}[b]{0.242\linewidth}
						\includetikz[width=\linewidth]{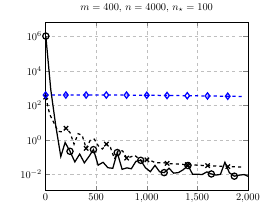}
						\includetikz[width=\linewidth]{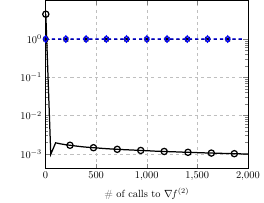}%
						\caption[]{\(\cost_1 = \|\cdot\|_1\)}
						\label{fig:linverseL1}%
					\end{subfigure}
					\caption[]{Linear inverse problems of \cref{sec:linverse} with minimum \(\ell^p\)-norm solution, \(p=1,2\).}
				\end{figure}

				In a series of experiments we consider the special cases of the following problem
				\begin{subequations}\label{eq:linverse}
					\begin{align}
						\minimize_{x\in\R^n}{} &~ \cost_1(x)
					\\
						\stt{} &~ x \in \argmin_{w\in\R^n}\tfrac12\|Aw -b\|^2,
					\end{align}
				\end{subequations}
				where \(A\in \R^{m\times n}\) and \(b\in \R^m\) are generated based on the procedure described in \cite[\S6]{nesterov2013gradient}, and \(n_\star\) denotes the number of nonzero elements of the solution.
				For the upper level cost \(\cost_1\), we consider two sets of experiments:
				\begin{enumerate}
				\item
					\(\cost_1 = \tfrac12\|{}\cdot{}\|^2\), corresponding to the Moore--Penrose solution, see \cref{fig:linverseL2};
				\item
					least \(\ell^1\)-norm solutions corresponding to \(\cost_1 = \|{}\cdot{}\|_1\), see \cref{fig:linverseL1} (the behavior of the algorithms is consistent with these plots for other values of \(m,n\)).
				\end{enumerate}
				As also done for the logistic regression problems, for \refadabim*, \refstabim*{} and \bisg{} we set \(\g_1 = \cost_1\) and \(\f_1 \equiv0\); for other methods the (smooth) upper level cost is captured by using \(\f_1\) with \(L_{\f_1} = 1\) and its strong convexity modulus equal to \(1\) (in the case of \bigsam{}).
				For methods that require Lipschitz modulus of \(\nabla \f_2\), \(L_{\f_2} = \|A\|^2\).

			\subsubsection{Solution of integral equations}
				We consider the solution of integral equations using the setting described in \cite[\S5.2]{beck2014first}.
				The corresponding bilevel problem is the following:
				\begin{subequations}\label{eq:IntEq}
					\begin{align}
						\minimize_{x\in\R^n}{} &~ \tfrac12\|x\|^2_Q
					\\
						\stt{} &~ x \in \argmin_{w\geq0}\tfrac12\|Aw -b\|^2.
					\end{align}
				\end{subequations}
				The data matrix \(A\) in \eqref{eq:linverse} is generated using \emph{philips, foxgood, baart} functions.
				Let \(L\) denote the discrete gradient operator, and let \(Q_1 = \trans LL\) and \(Q = Q_1 + \I\).
				In the simulations for \bigsam{} and \solodov, the upper level cost is captured using \(\f_1\), while for \refadabim*, \refstabim* and \bisg{} we used \(\f_1 = \tfrac12\langle x,Q_1x \rangle\) and \(\g_1 = \tfrac12\|{}\cdot{}\|^2\).
				\begin{figure}[h]
					\includetikz[width=\linewidth]{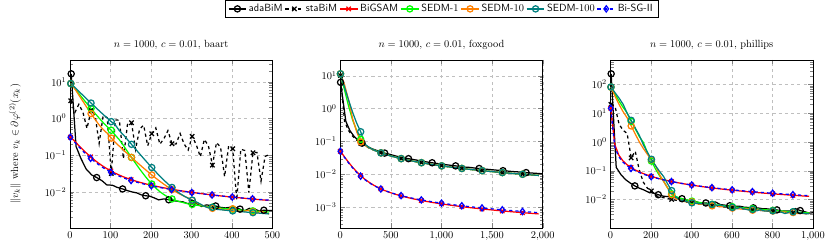}
					\includetikz[width=\linewidth]{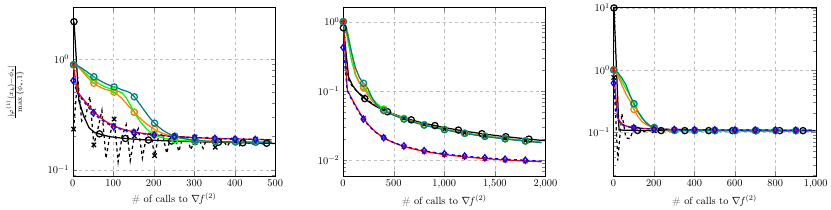}%
					\caption{Solution of integral equations}%
					\label{fig:inteq}%
				\end{figure}
				(As a rule of thumb, considering formulating the problem using the proximable term is  preferable, a tweak that only our methods and \bisg{} can take advantage of, being the only ones that allow proximable terms in the upper level.)
				By using the calculus rule of \cite[Prop. 24.8(i)]{bauschke2017convex}, the proximal mapping of \(\gk=\tfrac{\sigk}{2}\|{}\cdot{}\|^2+\g_2\) is given by
				\[
					\prox_{\alphk \gk}(u)
				=
					\prox_{\frac{\alphk}{1+\alphk\sigk}\g_2}\left(\tfrac{1}{1 + \alphk\sigk}u\right).
				\]
				Multiplications by \(Q_1\) involved in calls to \(\nabla\f_1\) can efficiently be handled through abstract linear operators and are ignored in the gradient calls count for all methods.
				In order to compare the methods in a fair manner, in addition to the deviation of the upper level cost from \(\costinf*\), we also plot a measure of optimality for the lower level.
				For example, in the case of \refadabim*, given that \(\g_1 = \tfrac12\|{}\cdot{}\|^2\), it is of immediate verification that
				\[
					v^k
				\coloneqq
					\tfrac{1}{\alphk*}(x_k-x_{k+1})-\sigk x_{k+1}+\nabla\f_2(x_{k+1})-\nabla\fk*(x_k)
				\in
					\partial \cost_2(x_{k+1}).
				\]
				Similar computation applies to the other methods that are included in the comparisons.


	\section{Conclusions}\label{sec:conclusions}

		This paper considered structured bilevel problems where both the upper and lower level minimizations are split as the sum of a nonsmooth and a (locally) Lipschitz differentiable function.
		A convergence recipe was developed for proximal gradient updates treating global and local Lipschitzian settings in a unified fashion. The aforementioned recipe depends on three properties for the stepsizes and involves a carefully designed adaptive scheme that builds upon and generalizes {\algnamefont adaPGM} \cite[Alg. 1]{latafat2023adaptive} to the bilevel setting.
		Notably, while in the locally Lipschitz setting our scheme involves a linesearch, it prescribes a suitable initialization for the linesearch based on Barzilai-Borwein type estimates, leading to much larger stepsizes compared to existing methods and considerably fewer backtracks in practice.
		Finally, the favorable convergence properties of the method were confirmed through a series of numerical simulations.
		Future research directions involve designing adaptive strategies for the inverse penalty parameters \(\sigk\),  deriving stopping criteria, and extensions to non-simple and possibly nonconvex bilevel settings.
		Relaxing the assumptions to (local) H\"older continuity of the gradients of the smooth terms is another promising direction which can leverage on recent developments on adaptive schemes \cite{oikonomidis2024adaptive}.


\ifarxiv
\else
	\begin{acknowledgements}
		\TheFunding
	\end{acknowledgements}
\fi

	\small
	\appendix
		\crefalias{section}{appendix}%
		\section{Appendix}\label{sec:appendix}

			\begin{appendixproof}{thm:stabim:convergence}
				Once \cref{cond:alphk,cond:alphk:nu,cond:alphamin} are verified, all the claims follow from \cref{thm:PG:barPhi,thm:recipe}.
				\Cref{cond:alphamin} is trivially satisfied for \refstabim* due to the underlying global Lipschitz continuity, and since \(\seq{\sigk}\) is decreasing.
				In what follows we consider the iterates generated by \refstabim*{} with \(\sigma_1=\sigma\) for some initial inverse penalty \(\sigma>0\). We show that for every \(k\) it holds that
				\[
					\alphk*\lk*\leq\nu
				\quad\text{and}\quad
					\alpha_1 = \tfrac{\nu}{\sigma_0L_{\f_1}+L_{\f_2}}\leq\alphk*\leq\min\set{\balphk*, \alpha_{\rm max}},
				\]
				where \(\balphk*\) is as in \eqref{eq:adabim:alphk*_0} with \(\alpha_{\rm max}\geq\frac{\nu}{L_{\f_2}}\) and initialization \(\alpha_{-1}=\alpha_0=\alpha_{\rm max}\) and \(\sigma_{-1}=\sigma_0=\sigma\).
				Since \(\lk*\leq\sigk* L_{\f_1}+L_{\f_2}\), the first inequality is obvious.
				Similarly, the second inequality follows from the fact that \(\seq{\sigk*}\) is decreasing, and thus so is \(\seq{\alphk*}\).
				Notice further that \(\alphk*\leq\frac{\nu}{L_{\f_2}}\leq\alpha_{\rm max}\); moreover, since
				\(\sigk\nabla\f_1+\nabla\f_2\) is globally Lipschitz with modulus \(\sigk L_{\f_1}+L_{\f_2}\), by cocoercivity it holds that
				\[
					\alphk\tfrac{\Lk^2}{\lk}
				=
					\alphk
					\tfrac{
						\|\nabla\fk(x_{k-1})-\nabla\fk(x_k)\|^2
					}{
						\innprod{\nabla\fk(x_{k-1})-\nabla\fk(x_k)}{x_{k-1}-x_k}
					}
				\leq
					\alphk(\sigk L_{\f_1}+L_{\f_2})
				=
					\nu
				<
					1,
				\]
				which as argued in \eqref{eq:ck} implies that the second term in \eqref{eq:adabim:alphk*_0} is infinite.
				Since, as already argued, \(\alphk\leq\alpha_{\rm max}\),
				to conclude it remains to show that \(\alphk*\) is also smaller than the first term in \eqref{eq:adabim:alphk*_0}.
				To this end, observe that
				\[
					\rhok*
				\defeq
					\frac{\sigk*\alphk*}{\sigk\alphk}
				=
					\frac{
						L_{\f_1}+\sigk^{-1}L_{\f_2}
					}{
						L_{\f_1}+\sigk*^{-1}L_{\f_2}
					}
				\leq
					1,
				\]
				where the inequality follows from the fact that \(\seq{\sigk*}\) is decreasing.
				Similarly,
				\(
					\rhok*
				=
					\tfrac{\sigk*\alphk*}{\sigk\alphk}
				\geq
					\tfrac{\sigk*}{\sigk}
				\geq
					\tfrac34
				\),
				where the first inequality follows from the fact that \(\seq{\alphk*}\) is increasing and the second one from the constraints on \(\sigk*\) prescribed at \cref{alg:stabim:sigk*}.
				Overall, it follows that \(\rhok*\in[\nicefrac34,1]\) holds for every \(k\).
				Therefore,
				\[
					\sqrt{
						\tfrac{\sigk}{\sigma_{k-1}}\bigl(1+\rhok\bigr)
					}
				\geq
					\sqrt{
						\tfrac{3}{4}\left(1+\tfrac{3}{4}\right)
					}
				>
					1
				\geq
					\rhok*
				\defeq
					\tfrac{\sigk*\alphk*}{\sigk\alphk},
				\]
				where the first inequality again follows from the bounds on \(\sigk*\) at \cref{alg:stabim:sigk*}.
				Rearranging yields the sought inequality
				\(
					\alphk*
				\leq
					\sqrt{\tfrac{\sigk}{\sigma_{k-1}}\bigl(1+\rhok\bigr)}
					\tfrac{\sigk}{\sigk*}
					\alphk
				\).
			\end{appendixproof}

			\begin{appendixproof}{thm:adabim:convergence}~
				We first state a simple lemma to justify the enlargement in the definition of the set \(\V\).

				\begin{lemma}\label{thm:r}%
					Let \(\func{h}{\R^n}{\Rinf}\) be proper lsc and convex, and given \(0<\alpha^+\leq\alpha^-\) and \(d\in\R^n\) let
					\(
						z^{\pm}\coloneqq\prox_{\alpha^{\pm}h}(x-\alpha^{\pm}d)
					\).
					Then, denoting \(\eta\coloneqq\nicefrac{\alpha^+}{\alpha^-}\in(0,1]\), it holds that
					\[
						\|z^--z^+\|
					\leq
						\tfrac{1-\eta}{\eta}
						\|x-z^+\|
					\quad\text{and}\quad
						\|z^--z^+\|^2
					\leq
						\tfrac{1-\eta}{1+\eta}
						\|x-z^-\|^2
						-
						\tfrac{1-\eta}{1+\eta}
						\|x-z^+\|^2.
					\]
				\end{lemma}
				\begin{proof}
					The proximal characterization of \(z^\pm\) reads
					\[
						\alpha^\pm h(z^\pm)
						+
						\tfrac{1}{2}\|z^\pm-x+\alpha^\pm d\|^2
					\leq
						\alpha^\pm h(y)
						+
						\tfrac{1}{2}\|y-x+\alpha^\pm d\|^2
						-
						\tfrac{1}{2}\|y-z^\pm\|^2
					\quad
						\forall y\in\R^n.
					\]
					By considering \(y=z^\mp\) and summing the resulting inequalities we obtain
					\begin{align*}
					&
						\tfrac{1}{2}\|z^+-x+\alpha^+ d\|^2
						+
						(\alpha^--\alpha^+)
						\bigl(h(z^-)-h(z^+)\bigr)
						+
						\tfrac{1}{2}\|z^--x+\alpha^- d\|^2
					\\
					\leq{} &
						\tfrac{1}{2}\|z^--x+\alpha^+ d\|^2
						+
						\tfrac{1}{2}\|z^+-x+\alpha^- d\|^2
						-
						\|z^+-z^-\|^2,
					\end{align*}
					which after expanding the squares and suitably rearranging results in
					\begin{align*}
						\tfrac{1}{\alpha^--\alpha^+}
						\|z^+-z^-\|^2
					\leq{} &
						h(z^+)-h(z^-)
						+
						\innprod{d}{z^+-z^-};
					\shortintertext{%
						using the fact that \(\tilde\nabla h(z^+)\coloneqq\frac{1}{\alpha^+}(x-z^+)-d\in\partial h(z^+)\) we may further upper bound this as
					}%
					\leq{} &
						\innprod{\tilde\nabla h(z^+)+d}{z^+-z^-}
					=
						\tfrac{1}{\alpha^+}
						\innprod{x-z^+}{z^+-z^-}.
					\end{align*}
					The first inequality in the statement now follows from the Cauchy-Schwarz inequality.
					The second one instead follows by expanding the inner product into three square norms and suitably rearranging.
				\end{proof}

				\begin{proofitemize}[topsep=0pt]
				\item \ref{thm:adabim:WD}~
					We begin by showing that the stepsize sequence is well defined and strictly positive.
					Let \(\alphk*_i\), \(x_{k+1,i}\), and \(\lk*_i\) respectively denote the value of \(\alphk*\), \(x_{k+1}\), and \(\lk*\) after \(i\) many backtracks, \(i\geq 0\);
					in particular,
					\(
						\alphk*_i
					=
						\min\set{\eta^i\balphk*,\alpha_{\rm max}}
					\)
					and
					\(
						x_{k+1,i}
					=
						\prox_{\alphk*_i\gk*}\left(x_k-\alphk*_i\nabla\fk*(x_k)\right)
					\).
					Let \(i_{k+1}\geq0\) denote the number of backtracks, or, equivalently, of failed attempts; we will show that \(i_{k+1}\) is finite for every \(k\); by construction, this will imply that \(x_{k+1}=x_{k+1,i_{k+1}}\).
					All the attempts \(\set{x_{k+1,i}}[\N\ni i\leq i_{k+1}]\) remain in a convex and compact set \(\V_{k+1}\)
					over which \(\nabla\fk*\) has finite Lipschitz modulus, be it \(L_{\fk*,\V_{k+1}}\);
					as such, one has that \(\alphk*_i\lk*_i\leq\alphk*_0\eta^iL_{\fk*,\V_{k+1}}\to0\) as \(i\to\infty\), implying that \eqref{eq:LS} is satisfied for \(i\) large enough.

					We thus proceed by induction to show that each iteration is well defined, that is, that \(\balphk*\) is (well defined and) strictly positive.
					Equivalently, in view of the update \eqref{eq:adabim:alphk*_0} it suffices to show that \(1-4\left(1-\tfrac{\sigk}{\sigma_{k-1}}\right)\alphk\lk^{(2)}>0\) holds for all \(k\).
					For \(k=0\) this is true because \(1-\frac{\sigma_0}{\sigma_{-1}}=0\) by initialization.
					Suppose that the claim holds for \(k\); then,
					\(
						\alphk\lk^{(2)}
					\leq
						\alphk\bigl(\sigk\lk^{(1)}+\lk^{(2)}\bigr)
					=
						\alphk\lk
					\leq
						\nu
					\),
					where the last inequality owes to the linesearch condition \eqref{eq:LS} at the previous step, which holds by inductive hypothesis.
					Therefore,
					\begin{equation}\label{eq:sqrt>0}
						1
						-
						4\left(1-\tfrac{\sigk}{\sigma_{k-1}}\right)
						\alphk\lk^{(2)}
					\geq
						1
						-
						4\nu\left(1-\tfrac{\sigk}{\sigma_{k-1}}\right)
					\geq
						1 - \nu
					>
						0,
					\end{equation}
					where the last inequality owes to the bound \(\sigk\geq\frac{3}{4}\sigma_{k-1}\) prescribed at \cref{state:adabim:alphk*_0}.

					This concludes the proof of the well definedness of the iterations.
					Notice also that \cref{cond:alphk,cond:alphk:nu} hold by construction, thereby ensuring through \cref{thm:PG:boundedness,thm:PG:barPhi} that \(\seq{x^k}\) is bounded and \(\|x^{k+1}-x^k\|\to0\).
					In particular, \(r=\sup_{k\in\N}\|x^{k+1}-x^k\|<\infty\), implying that \(\V\) as in the statement is bounded.
					The enlargement in its definition also guarantees that, in addition to all the iterates \(x_k\), \(\V\) also contains \(x_{k,i_k^-}\) where \(i_k^-=[i_k-1]_+\) (\(\V\) contains every \(x_k\) and all the last failed attempts whenever the linsearch is not passed at the first trial).
					This follows from the first inequality in \cref{thm:r} with \(x=x_k\), \(h=\gk*\), \(d=\nabla\fk*(x_k)\), \(\alpha^+=\alphk*\) and \(\alpha^-=\alphk*_{i_k^-}=\alphk*/\eta\), for which \(z^+=x_{k+1}\) and \(z^-=x_{k+1,i_{k+1}^-}\).

					Whenever \(\balphk*< \tfrac{\nu}{L_{f_0,\V}}\) the initial stepsize \(\balphk*\) already complies with \eqref{eq:LS}, and therefore \(\alphk*=\min\set{\balphk*,\alpha_{\rm max}}=\balphk*\) (the last identity follows from the fact that \(\alpha_{\rm max}\geq\frac{1}{\ell_0}\geq\frac{1}{L_{f_0,\V}}\)).
					Otherwise, the failure of the second-last backtrack implies that
					\(
						\alphk*
					=
						\eta\alphk*_{i_k-1}
					>
						\tfrac{\eta\nu}{L_{f_0,\V}}
					\).
					Either way, this shows that
					\begin{equation}\label{eq:ahata}
						\alphk*
					\geq
						\min\set{
							\balphk*,
							\tfrac{\eta\nu}{L_{f_0,\V}}
						}
					\quad
						\forall k.
					\end{equation}
					We show by induction that \(\alphk\geq \alpha_{\rm min}\).
					For \(k=1\), by the choice \(\alpha_{-1}\) used in the initialization of \refadabim*{} and since \(\sigma_0 = \sigma_{-1}\), we have
					\begin{align*}
						\widehat\alpha_{1}
					={} &
						\tfrac{\sigma_{0}}{\sigma_{1}}
						\alpha_0
						\min\set{
							\sqrt{1+\tfrac{\alpha_0}{\alpha_{-1}}}
							,\,
							\tfrac{
								1
							}{
								2\sqrt{[\alpha_0^2 L_0^2 - \alpha_0\ell_0]_+}
							}
						}
					\geq
						\min\set{
							\alpha_0\sqrt{1+\tfrac{\alpha_0}{\alpha_{-1}}}
							,\,
							\tfrac{
								1
							}{
								2L_0
							}
						}
					\\
					\geq{} &
						\left.
							\begin{ifcases}
								\min\set{
									\sqrt{2}\alpha_0
									,\,
									\tfrac{
										1
									}{
										2L_0
									}
								}
							&
								\alpha_0\ell_0\geq\nicefrac{1}{2}
							\\
								\min\set{
									\tfrac1{\ell_0}
									,\,
									\tfrac{
										1
									}{
										2L_0
									}
								}
							\otherwise
							\end{ifcases}
						\right\}
					\geq
						\tfrac1{2L_{f_0, \V}}
					\geq
						\alpha_{\rm min}.
					\numberthis\label{eq:alpha_1}
					\end{align*}
					Therefore, by \eqref{eq:ahata} also \(\alpha_1 \geq \alpha_{\rm min}\).
					Suppose that the claim holds up to \(k\). We consider two cases.
					\begin{itemize}[leftmargin = *, topsep=0pt,itemsep=0pt]
					\item
						If \(\balphk*\) equals the second element in \eqref{eq:adabim:alphk*_0}, then
						\[
							\balphk*
						=
							\alphk
								\tfrac{\sigk}{\sigk*}
							\tfrac{
								\sqrt{
									1
									-
									4\left(1-\frac{\sigk}{\sigma_{k-1}}\right)
									\alphk\lk^{(2)}
								}
							}{
								2\sqrt{\alphk^2 \Lk^2 - \alphk\lk}
							}
						\geq
							\alphk
							\tfrac{
								\sqrt{1-\nu}
							}{
								2\sqrt{\alphk^2 \Lk^2 - \alphk\lk}
							}
						\geq
							\tfrac{
								\sqrt{1-\nu}
							}{
								2L_k
							}
						\geq
							\tfrac{
								\sqrt{1-\nu}
							}{
								2L_{f_0,\V}
							}
						\]
						which combined with \eqref{eq:ahata} yields that
						\(
							\alphk*
						\geq
							\min\set{
								\tfrac{
									\sqrt{1-\nu}
								}{
									2L_{f_0,\V}
								}
							,
								\frac{\eta\nu}{L_{f_0,\V}}
							}
							\geq \alpha_{\rm min}
						\).

					\item
						Suppose that \(\balphk*\) equals the first element in \eqref{eq:adabim:alphk*_0}.
						If \(\alphk* < \balphk*\), then by \eqref{eq:ahata} the lower bound \(\alphk* \geq \tfrac{\eta\nu}{L_{f_0,\V}}\geq \alpha_{\rm min}\) holds.
						If instead \(\alphk* = \balphk*\), then
						\begin{equation}\label{eq:alphk*=balphk*}
							\alphk*
						=
							\balphk*
						=
							\tfrac{\sigk}{\sigk*}
							\alphk
							\sqrt{\tfrac{\sigk}{\sigma_{k-1}}\bigl(1+\rhok\bigr)}
						\geq
							\alphk
							\sqrt{\tfrac{3}{4}\bigl(1+\rhok\bigr)}.
						\end{equation}

						We consider two subcases.
						\begin{itemize}[topsep=0pt,itemsep=0pt]
						\item
							If \(\alphk<\balphk\), then, as argued above, \(\alphk \geq \tfrac{\eta\nu}{L_{f_0,\V}}\) and by \eqref{eq:alphk*=balphk*}
							\(
								\alphk*
							\geq
								\alphk
								\sqrt{\tfrac{3}{4}\bigl(1+\rhok\bigr)}
							\geq
								\tfrac{\sqrt{3}}{2}
								\tfrac{\eta\nu}{L_{f_0,\V}}
							\geq
								\alpha_{\rm min}
							\).

						\item
							If instead \(\alphk = \balphk\), then observe that
							\(
								\rhok
							\defeq
								\frac{\sigk\alphk}{\sigma_{k-1}\alpha_{k-1}}
							=
								\frac{\sigk\balphk}{\sigma_{k-1}\alpha_{k-1}}
							=
								\sqrt{\tfrac{\sigma_{k-1}}{\sigma_{k-2}}\bigl(1+\rho_{k-1}\bigr)}
							\),
							hence
							\begin{align*}
								\alphk* = \balphk*
							={} &
								\tfrac{\sigk}{\sigk*}
								\alphk
								\sqrt{\tfrac{\sigk}{\sigma_{k-1}}\bigl(1+\rhok\bigr)}
							\\
							={} &
								\tfrac{\sigk}{\sigk*}
								\alphk
								\sqrt{
									\tfrac{\sigk}{\sigma_{k-1}}
									\Bigl(
										1
										+
										\sqrt{\tfrac{\sigma_{k-1}}{\sigma_{k-2}}\bigl(1+\rho_{k-1}\bigr)}
									\Bigr)
								}
							\geq
								\alphk
								\sqrt{
									\tfrac34
									\bigl(
										1
										+
										\tfrac{\sqrt{3}}{2}
									\bigr)
								}
							>
								\alphk
							\geq
								\alpha_{\rm min},
							\end{align*}
							where the first inequality uses
							\(\sigma_j\geq\sigma_{j+1}\geq\frac{3}{4}\sigma_j\) for every \(j\),
							and the last one holds by induction.
						\end{itemize}
						We showed that in either case \(\alphk* \geq \alpha_{\rm min}\) as claimed.
					\end{itemize}

				\item \ref{thm:adabim:sublinear} and \ref{thm:adabim:conv}~
					Having established \cref{cond:alphamin} and since \cref{cond:alphk,cond:alphk:nu} hold by design, both assertions follow immediately from \cref{thm:PG:barPhi,thm:recipe}.
				\qedhere
				\end{proofitemize}
			\end{appendixproof}

	\phantomsection
	\addcontentsline{toc}{section}{References}%
	\bibliographystyle{plain}
	\bibliography{Bibliography_abbr.bib}

\begin{thebibliography}{10}

\bibitem{attouch1996viscosity}
H.~Attouch.
\newblock Viscosity solutions of minimization problems.
\newblock {\em SIAM J. Optim.}, 6(3):769--806, 1996.

\bibitem{bahraoui1994convergence}
M.A. Bahraoui and B.~Lemaire.
\newblock Convergence of diagonally stationary sequences in convex
  optimization.
\newblock {\em Set-Valued Anal.}, 2:49--61, 1994.

\bibitem{barzilai1988two}
J.~Barzilai and J.M. Borwein.
\newblock Two-point step size gradient methods.
\newblock {\em IMA J. Numer. Anal.}, 8(1):141--148, 1988.

\bibitem{bauschke2017convex}
H.H. Bauschke and P.L. Combettes.
\newblock {\em Convex analysis and monotone operator theory in {H}ilbert
  spaces}.
\newblock CMS Books Math. Springer, 2017.

\bibitem{beck2017first}
A.~Beck.
\newblock {\em First-Order Methods in Optimization}.
\newblock SIAM, Philadelphia, PA, 2017.

\bibitem{beck2014first}
A.~Beck and S.~Sabach.
\newblock A first order method for finding minimal norm-like solutions of
  convex optimization problems.
\newblock {\em Math. Program.}, 147(1-2):25--46, 2014.

\bibitem{bigi2022combining}
G.~Bigi, L.~Lampariello, and S.~Sagratella.
\newblock Combining approximation and exact penalty in hierarchical
  programming.
\newblock {\em Optim.}, 71(8):2403--2419, 2022.

\bibitem{borsos2020coresets}
Z.~Borsos, M.~Mutny, and A.~Krause.
\newblock Coresets via bilevel optimization for continual learning and
  streaming.
\newblock {\em Adv. Neural Inf. Process. Syst.}, 33:14879--14890, 2020.

\bibitem{cabot2005proximal}
A.~Cabot.
\newblock Proximal point algorithm controlled by a slowly vanishing term:
  applications to hierarchical minimization.
\newblock {\em SIAM J. Optim.}, 15(2):555--572, 2005.

\bibitem{chang2011libsvm}
C.~Chang and C.~Lin.
\newblock {LIBSVM}: A library for support vector machines.
\newblock {\em ACM Trans. Intell. Syst. Technol. (TIST))}, 2:1--27, 2011.

\bibitem{dempe2002foundations}
S.~Dempe.
\newblock {\em Foundations of bilevel programming}.
\newblock Springer Science \& Business Media, 2002.

\bibitem{dempe2020bilevel}
S.~Dempe.
\newblock Bilevel optimization: theory, algorithms, applications and a
  bibliography.
\newblock {\em Bilevel Optim.: Adv. Next Challenges}, pages 581--672, 2020.

\bibitem{doron2022methodology}
L.~Doron and S.~Shtern.
\newblock Methodology and first-order algorithms for solving nonsmooth and
  non-strongly convex bilevel optimization problems.
\newblock {\em Math. Program.}, pages 1--38, 2022.

\bibitem{facchinei2014vi}
F.~Facchinei, J.~Pang, G.~Scutari, and L.~Lampariello.
\newblock {VI}-constrained hemivariational inequalities: distributed algorithms
  and power control in ad-hoc networks.
\newblock {\em Math. Program.}, 145(1-2):59--96, 2014.

\bibitem{franceschi2018bilevel}
L.~Franceschi, P.~Frasconi, S.~Salzo, R.~Grazzi, and M.~Pontil.
\newblock Bilevel programming for hyperparameter optimization and
  meta-learning.
\newblock In {\em Int. Conf. Mach. Learn.}, pages 1568--1577. PMLR, 2018.

\bibitem{garrigos2018iterative}
G.~Garrigos, L.~Rosasco, and S.~Villa.
\newblock Iterative regularization via dual diagonal descent.
\newblock {\em J. Math. Imaging Vision}, 60:189--215, 2018.

\bibitem{grazzi2023bilevel}
R.~Grazzi, M.~Pontil, and S.~Salzo.
\newblock Bilevel optimization with a lower-level contraction: Optimal sample
  complexity without warm-start.
\newblock {\em J. Mach. Learn. Res.}, 24(167):1--37, 2023.

\bibitem{guan2023first}
W.~Guan and W.~Song.
\newblock A first-order method for solving bilevel convex optimization problems
  in {B}anach space.
\newblock {\em Optim.}, pages 1--26, 2023.

\bibitem{helou2017subgradient}
E.S. Helou and L.E. Sim{\~o}es.
\newblock $\epsilon$-subgradient algorithms for bilevel convex optimization.
\newblock {\em Inverse Probl.}, 33(5):055020, apr 2017.

\bibitem{hong2023two}
M.~Hong, H.~Wai, Z.~Wang, and Z.~Yang.
\newblock A two-timescale stochastic algorithm framework for bilevel
  optimization: Complexity analysis and application to actor-critic.
\newblock {\em SIAM J. Optim.}, 33(1):147--180, 2023.

\bibitem{jiang2023conditional}
R.~Jiang, N.~Abolfazli, A.~Mokhtari, and E.Y. Hamedani.
\newblock A conditional gradient-based method for simple bilevel optimization
  with convex lower-level problem.
\newblock In {\em Int. Conf. Artif. Intell. Stat.}, pages 10305--10323. PMLR,
  2023.

\bibitem{kaushik2021method}
H.D. Kaushik and F.~Yousefian.
\newblock A method with convergence rates for optimization problems with
  variational inequality constraints.
\newblock {\em SIAM J. Optim.}, 31(3):2171--2198, 2021.

\bibitem{lampariello2022solution}
L.~Lampariello, G.~Priori, and S.~Sagratella.
\newblock On the solution of monotone nested variational inequalities.
\newblock {\em Math. Methods Oper. Res.}, 96(3):421--446, 2022.

\bibitem{latafat2023convergence}
P.~Latafat, A.~Themelis, and P.~Patrinos.
\newblock On the convergence of adaptive first order methods: proximal gradient
  and alternating minimization algorithms.
\newblock {\em \arXivLink{2311.18431}}, 2023.

\bibitem{latafat2023adaptive}
P.~Latafat, A.~Themelis, L.~Stella, and P.~Patrinos.
\newblock Adaptive proximal algorithms for convex optimization under local
  {L}ipschitz continuity of the gradient.
\newblock {\em \arXivLink{2301.04431}}, 2023.

\bibitem{malitsky2020adaptive}
Y.~Malitsky and K.~Mishchenko.
\newblock Adaptive gradient descent without descent.
\newblock In {\em Proc. 37th Int. Conf. Mach. Learn.}, volume 119, pages
  6702--6712. PMLR, 13- 2020.

\bibitem{merchav2023convex}
R.~Merchav and S.~Sabach.
\newblock Convex bi-level optimization problems with non-smooth outer objective
  function, 2023.

\bibitem{moudafi2000viscosity}
A.~Moudafi.
\newblock Viscosity approximation methods for fixed-points problems.
\newblock {\em J. Math. Anal. Appl.}, 241(1):46--55, 2000.

\bibitem{nesterov2013gradient}
Y.~Nesterov.
\newblock Gradient methods for minimizing composite functions.
\newblock {\em Math. Program.}, 140(1):125--161, aug 2013.

\bibitem{oikonomidis2024adaptive}
K.A. Oikonomidis, E.~Laude, P.~Latafat, A.~Themelis, and P.~Patrinos.
\newblock Adaptive proximal gradient methods are universal without
  approximation.
\newblock {\em \arXivLink{2402.06271}}, 2024.

\bibitem{pedregosa2016hyperparameter}
F.~Pedregosa.
\newblock Hyperparameter optimization with approximate gradient.
\newblock In {\em Int. Conf. Mach. Learn.}, pages 737--746. PMLR, 2016.

\bibitem{peypouquet2012coupling}
J.~Peypouquet.
\newblock Coupling the gradient method with a general exterior penalization
  scheme for convex minimization.
\newblock {\em J. Optim. Theory Appl.}, 153:123--138, 2012.

\bibitem{rajeswaran2019meta}
A.~Rajeswaran, C.~Finn, S.M. Kakade, and S.~Levine.
\newblock Meta-learning with implicit gradients.
\newblock {\em Adv. neural Inf. Process. Syst.}, 32, 2019.

\bibitem{rockafellar2011variational}
R.T. Rockafellar and R.J. Wets.
\newblock {\em Variational analysis}, volume 317.
\newblock Springer, 2011.

\bibitem{sabach2017first}
S.~Sabach and S.~Shtern.
\newblock A first order method for solving convex bilevel optimization
  problems.
\newblock {\em SIAM J. Optim.}, 27(2):640--660, 2017.

\bibitem{solodov2007bundle}
M.V. Solodov.
\newblock A bundle method for a class of bilevel nonsmooth convex minimization
  problems.
\newblock {\em SIAM J. Optim.}, 18(1):242--259, 2007.

\bibitem{solodov2007explicit}
M.V. Solodov.
\newblock An explicit descent method for bilevel convex optimization.
\newblock {\em J. Convex Anal.}, 14(2):227, 2007.

\bibitem{themelis2019acceleration}
A.~Themelis, M.~Ahookhosh, and P.~Patrinos.
\newblock On the acceleration of forward-backward splitting via an inexact
  {N}ewton method.
\newblock In H.H. Bauschke, R.S. Burachik, and D.R. Luke, editors, {\em
  Splitting Algorithms, Modern Operator Theory, and Applications}, pages
  363--412. Springer, 2019.

\bibitem{xu2004viscosity}
H.~Xu.
\newblock Viscosity approximation methods for nonexpansive mappings.
\newblock {\em J. Math. Anal. Appl.}, 298(1):279--291, 2004.

\end{thebibliography}

\end{document}